\input amstex
\documentstyle{amsppt}

\vcorrection{-.4in}\hcorrection{.2 in}

\input xypic
\input epsf

\define\CDdashright#1#2{&\,\mathop{\dashrightarrow}\limits^{#1}_{#2}\,&}
\define\CDdashleft#1#2{&\,\mathop{\dashleftarrow}\limits^{#1}_{#2}\,&}

\def\P{\Bbb P}

\def\C{\Bbb C}

\def\Q{\Bbb Q}

\def\Til#1{{\widetilde{#1}}}

\def\PGL{\text{PGL}}

\def\im{\text{\rm im}\,}

\def\siltable#1.{
\vbox{\tabskip=0pt \offinterlineskip
\halign to360pt{\strut##& ##\tabskip=1em plus2em&
\hfil##\hfil& \vrule##&
\hfil##\hfil& \vrule##&
\hfil##\hfil& \vrule##&
\hfil##\hfil& \vrule\thinspace\vrule##&
\hfil##\hfil& \vrule\thinspace\vrule##&
\hfil##\hfil& ##\tabskip=0pt\cr
#1}}}


\magnification=1200
\CompileMatrices

\topmatter
\title Linear orbits of arbitrary plane curves\endtitle
\author Paolo Aluffi and Carel Faber\endauthor
\date December 1999\enddate
\address Mathematics Department, Florida State University,
Tallahassee FL 32306, U.S.A. \endaddress
\email aluffi\@math.fsu.edu \endemail
\address
Dept.~of Mathematics,
Oklahoma State University,
Stillwater OK 74078, U.S.A.
\endaddress
\email cffaber\@math.okstate.edu\endemail
\address
Department of Mathematics,
KTH,
100 44 Stockholm, Sweden
\endaddress
\email carel\@math.kth.se\endemail

\abstract The `linear orbit' of a plane curve of degree $d$ is its
orbit in $\P^{d(d+3)/2}$ under the natural action of $\PGL(3)$. In
this paper we obtain an algorithm computing the degree of the closure 
of the linear orbit of an arbitrary plane curve, and give explicit
formulas for plane curves with irreducible singularities. The main
tool is an intersection@-theoretic study of the projective normal cone
of a scheme determined by the curve in the projective space $\P^8$ of
$3\times 3$ matrices; this expresses the degree of the orbit closure
in terms of the degrees of suitable loci related to the limits of the
curve. These limits, and the degrees of the corresponding loci, have
been established in previous work.
\endabstract


\leftheadtext{Paolo Aluffi and Carel Faber}

\endtopmatter

\document
\footnote[]{Both authors gratefully acknowledge partial NSF support,
under grants DMS-9500843 and DMS-9801257.}


\head \S0. Introduction\endhead

The Gromov@-Witten invariants of $\P^2$ compute, roughly speaking, the
number of plane curves of given degree $d$ and genus $g$ containing the
appropriate number of general points. In recent years it has been
discovered that these invariants are coherently linked together by the
apparatus of quantum cohomology, which exposes their
structure as $d$ and $g$ are allowed to vary.

For {\it nonsingular\/} plane curves, however, these invariants do not
carry much information: the set of nonsingular curves of a given
degree $d$ is an open set of a projective space $\P^{d(d+3)/2}$, so
the corresponding invariant is simply~1. We can consider a more
refined question by fixing, as well as the degree $d$ (and therefore
the genus $g=\frac{(d-1)(d-2)}2$), the moduli class in $\Cal M_g$ of
the curve. What data determines then the corresponding invariant? Can
this invariant be effectively computed? Can other enumerative
invariants be computed for the set of nonsingular curves of given
degree and moduli class, such as the number of curves tangent to the
appropriate number of general lines?

In this article we fully answer these questions, and a natural
generalization of these questions to {\it arbitrary\/} (i.e., possibly
singular, reducible, nonreduced) plane curves of any degree.
The group $\PGL(3)$ of projective linear transformations of $\P^2$
acts naturally on the space $\P^{d(d+3)/2}$ parametrizing plane curves
of degree~$d$. Our main result is the computation of the degree of the
closure  in this space
of the orbit of an arbitrary plane curve (in char.~0).
The orbit closure of a curve is a natural object of study, and
its degree has a simple enumerative meaning: for a reduced curve
with finite stabilizer,
it counts the number of translates of
the curve which contain 8~given general points. For a nonsingular
curve, this is the invariant mentioned above. In this sense,
therefore, this problem is an isotrivial version of the problem of
computing Gromov@-Witten invariants.
Somewhat surprisingly, the enumerative geometers and the invariant
theorists of the 19th century
do not seem to have worked on this question.

The computation in this paper relies on our previous work on the
subject, where we have 
dealt with special curves: nonsingular curves
were in fact already treated in \cite{A-F2}; plane curves whose orbit
has dimension less than $\dim\PGL(3)=8$ are classified and studied in
\cite{A-F3}, \cite{A-F4}. We have also determined in \cite{A-F5} the
{\it limits\/}  of an arbitrary plane curve; these are the curves
appearing in the boundary of the orbit, that is, the
complement of the orbit in its closure. In the terminology 
of \cite{H-M} (p.~138) this solves the `isotrivial flat completion
problem' for plane curves. 

Our previous enumerative computations
relied on the explicit construction (by means of a sequence of
blow@-ups over the $\P^8$ of $3\times 3$ matrices) of smooth varieties
dominating the orbit closures. The case of an arbitrary curve appears
to be too complex for that approach, and we turn in this paper to a
more direct study of the projective normal cone of the base locus
(scheme) of the rational map
$$\P^8 \dashrightarrow \P^{d(d+3)/2}$$
extending the map $\PGL(3) @>>> \P^{d(d+3)/2}$ which surjects onto the
orbit of a given curve. Our study of limits of curves in \cite{A-F5}
allows us to express the degree of the orbit closure of a curve in
terms of enumerative information concerning curves in the boundary of the
orbit, also
available from our previous computations.

For an arbitrary curve, this provides us implicitly with an algorithm
computing the degree of the orbit closure. We illustrate this
algorithm in \S4 and \S5 on specific classes of curves. For example, a
surprisingly simple formula can be obtained to compute the effect on
the degree due to an {\it irreducible\/} singularity $p$ of a curve
(see Theorem~5.1)
in terms of the multiplicity of the curve at $p$, the order of contact with
the tangent line to the branch at $p$, and the Puiseux pairs
describing the singularity.

Of course many questions remain concerning orbit closures, for example
regarding their singularities 
(which curves have smooth orbit closure?~smooth
orbit closures of {\it configurations of points in $\P^1$\/} are
classified in \cite{A-F1}),
or other invariants such as Euler characteristic, Poincar\'e
polynomials, behavior in positive characteristic, etc.
\vskip 6pt

{\bf Acknowledgement.\/} 
It is a pleasure to thank Bill Fulton for several invitations to visit
the University of Chicago, where much of the work on this project was done.


\head \S1. The problem, and the approach\endhead

Let $C$ be a curve of degree $d$ in the projective plane $\P^2$ over
an algebraically closed field of characteristic~0; we may think of
$C$ as a point in the projective space $\P^N=\P(H^0(\P^2,\Cal O(d)))$,
where $N=d(d+3)/2$. The standard action of $\PGL(3)$ on $\P^2$ induces
a right action on $\P^N$; specifically, for $\varphi\in\PGL(3)$ we can
consider the {\it translate\/} of $C$ by $\varphi$: if $C$ has
equation $F(x_0:x_1:x_2)=0$, then its translate $C\circ\varphi$ has
equation
$$F(\varphi(x_0:x_1:x_2))=0\quad.$$
The action $\varphi\mapsto C\circ\varphi$ defines a map
$$c:\PGL(3) @>>> \P^N$$
whose image is what we call the {\it linear orbit\/} of $C$. Our aim
is the computation of the degree of the closure of this orbit, for an
arbitrary plane curve $C$, in terms of a description of the
irreducible components and the singularities of $C$. 

Our general approach is based on compactifying $\PGL(3)$ to the space
$\P^8$ of $3\times 3$ matrices, and considering the rational map
$$\P^8 \dashrightarrow \P^N$$
determined by $c$. If $\tilde c: \Til V @>>> \P^N$ is a map
resolving the indeterminacies of this rational map, so that the diagram
$$\diagram
{\Til V} \dto_{\pi} \drto^{\tilde c}\\
{\P^8} \xdashed[0,1]^{c}|>\tip & {\P^N} \\
\enddiagram$$
commutes, then the orbit closure of $C$ is the image of $\tilde c$. In
special but important cases one can in fact construct and study a
{\it nonsingular\/} such variety $\Til V$, by a suitable sequence of
blow@-ups along smooth centers over $\P^N$; this is carried out in
\cite{A-F2}, \cite{A-F3}, \cite{A-F4}. The work involved in the
construction of an explicit resolution of the orbit closure pays off in
terms of a simpler intersection@-theoretic set@-up, and 
opens the door to a more thorough study of the orbit closure.

Such a construction is however not available for an {\it arbitrary\/}
plane curve $C$. This is an indication of the fact that singularities
of a plane curve can be extremely complicated, and that the
orbit closure is highly sensitive to the local features of a curve. To
treat the general case, we resort then essentially to using the most
simple@-minded (but highly singular) variety $\Til V$ as above---we will
let $\Til V$ be the blow@-up of $\P^8$ along the base scheme $S$ of the
rational map $c$---and pay the price of a more complicated
intersection@-theoretic set@-up and of a careful local study of
degenerations of $C$. In the end we will be able to express the degree
of the orbit closure of $C$ in terms of enumerative information
concerning its {\it limits\/}, that is, the curves obtained as limits of
translates $C\circ\varphi$ as $\varphi$ approaches 
the base locus of $c$.
This enumerative information has been obtained in our previous
work; it relies on the explicit
resolution of the orbit closure of the limits.

In this section we describe our degeneration technique, and the
intersection theory formula we will use in the main computation.
The degree of the orbit closure is the intersection number
$$h^{\dim\tilde c(\Til V)}\cdot [\tilde c(\Til V)]\quad,$$
where $h$ denotes the hyperplane class in $\P^N$. Pulling back to
$\Til V$, we are then led to consider the class
$$h^{\dim\tilde c(\Til V)}\cap [\Til V]$$
(following common practice, we omit evident pull@-back notations);
in fact, in order not to fix from the start the dimension of the orbit
of $C$, we consider the class
$$\frac{[\Til V]}{c(\Cal O(-h))} = (1+h+h^2+h^3+\dots)\cap [\Til V]
\quad,$$
and its push@-forward to $\P^8$:
$$\pi_*\frac{[\Til V]}{c(\Cal O(-h))}=(1+a_1 H+ a_2 H^2+\dots)
\cap[\P^8]\quad,$$
where $H$ is the hyperplane class in $\P^8$, and $a_i$ is the degree
of $\pi_*(h^i\cap [\Til V])$. It is clear that
$$a_i=0\quad\text{for $i>\dim\tilde c(\Til V)$}$$
and that $a_{\dim\tilde c(\Til V)}$ equals the degree of the orbit
closure times the degree of the closure of the stabilizer of $C$ in
$\P^8$. We call this number the `predegree' of the orbit closure of
$C$, and the whole class written above, which we think of as a
polynomial in $H$, the {\it predegree polynomial\/} of (the orbit
closure of) $C$. 

\remark{Note} The `polynomials' appearing in this paper are therefore
nothing but classes in the Chow ring of $\P^8$. It will in fact be
convenient to take rational coefficients, so that our polynomials will
live in the ring $\Q[H]/(H^9)$. When manipulating polynomials we will
implicitly work in this ring; in particular, all operations are
truncated to $H^8$. This allows us some convenient abuse of language;
for example, 
$$\exp(dH)=1+d H+\frac{(dH)^2}2+\frac{(dH)^3}{3!}+\dots+
\frac{(dH)^8}{8!}$$
with our conventions.\endremark

Our objective then becomes the following: {\it compute the predegree
polynomial of an arbitrary plane curve $C$.} The degree of the orbit
closure of a curve $C$ is recovered from its predegree polynomial by
dividing the top nonzero coefficient by the degree of the closure of
the stabilizer of $C$. Predegree polynomials are a more natural object
of study, since they carry enumerative information independently of
the dimension of the orbit closure. 
The information in the predegree
polynomial is equivalent to the information in what we call the
{\it adjusted predegree polynomial (a.p.p.)}
$$\pi_*(ch(\Cal O(h))\cap [\Til V])=1+a_1 H+a_2\frac{H^2}2+ a_3
\frac{H^3}{3!}+\dots\quad.$$
Computing {\it adjusted\/} predegree polynomials often leads to
simpler formulas, so we focus on them in this paper. Adjusted
predegree polynomials for curves with small orbits
(i.e., of dimension $<8$) are computed in
\cite{A-F3}, \cite{A-F4}.

We can analyze the situation in a more general context. Let $V$ be any
variety, $\Cal L$ a line bundle on $V$, and $\Cal E\subset H^0(V, \Cal
L)$ a nonzero linear system.
These choices determine a rational map
$$\alpha: V \dashrightarrow \P^N=\P(\Cal E^\vee)\quad.$$
Let $S$ be the scheme@-theoretic intersection of the sections in $\Cal
E$, so that the base locus of $\alpha$ is the support of $S$, and (the
closure of) the graph $\Gamma$ of $\alpha$ can be identified with the
blow@-up $\Til V$ of $V$ along $S$. We let $E$ be the exceptional
divisor of the blow@-up, that is, the part of the graph over $S$:
$$\diagram
{E=\pi^{-1}S} \drto \rto & {\Gamma=\Til V}\drto^{\pi} \rrto && {V\times
\P^N}\dlto \drto\\ 
& S \rto & V \xdashed[0,2]^{\alpha}|>\tip && {\P^N}
\enddiagram$$
In other words, $E$ is a realization of the projective normal cone of
$S$ in $V$. Let now $\Til{\Cal L}$ denote the pull@-back to $\Gamma$ of
the hyperplane class in $\P^N$, and notice that if $\Cal E$ is
base@-point@-free to begin with (so $S=\emptyset$), then $\Til{\Cal
L}=\Cal L$ and the quantity corresponding to the adjusted predegree
polynomial is simply
$$\pi_*(ch(\Til{\Cal L})\cap [\Til V])=ch(\Cal L)\cap [V]
\quad\text{in $(A_*V)_\Q$.}\tag*$$
The following proposition
shows how to modify the fundamental class of $[V]$ in this formula to
account for the base locus $S$ of $\alpha$. The correction term will be
obtained from the cycle of~$E$:
$$[E]=m_1 [E_1]+\dots+m_r [E_r]\quad,$$
as follows. We denote by $h$ the hyperplane class in $\P^N$
and its pull@-backs (for example, $h=c_1(\Til{\Cal L})$ on
$\Gamma$); write $\ell=c_1(\Cal L)$, and let
$$L_i=\sum_{k\ge 0}\frac 1{k+1}\sum_{j=0}^k\frac{(-\ell)^{k-j}}
{j!(k-j)!}\pi_*(h^j\cap [E_i])$$
(so {\it a priori\/} the $L_i$ might have nonzero terms in all
dimensions from 0 to $\dim V-1$). Here is the main observation in this
section:
\proclaim{Proposition 1.1}
$$\pi_*(ch(\Til{\Cal L})\cap [\Til V])=ch(\Cal L)\cap\big([V]-(m_1
L_1+ \dots +m_r L_r)\big)\quad\text{in $(A_*V)_\Q$.}$$
\endproclaim
\demo{Proof} Note that $h=c_1(\Til{\Cal L})=\ell-e$, where $e$ is the
class of $E$ and as usual we omit obvious pull@-back notations.
Therefore
$$\align
\pi_* (ch(\Til{\Cal L})\cap [\Til V]) &=\pi_*(\exp(\ell-e)\cap [\Til
V]) =\exp(\ell)\cap \pi_*(\exp(-e)\cap[\Til V])\\
& =\exp(\ell)\cap ([V]-\pi_*(1-\exp(-e))\cap[\Til V])
\endalign$$
giving the correction term to the fundamental class as
$$\gather
-\pi_*(1-\exp(-e))\cap[\Til V]=-\pi_*\sum_{i\ge 0} \frac{(-e)^i} {(i+1)!}
\cap[E]\quad,\text{ that is}\\
-\pi_*\sum_{i\ge 0} \frac{(h-\ell)^i} {(i+1)!} \cap(m_1[E_1] +\dots
+m_r[E_r])\quad.
\endgather$$
The statement follows by expanding this expression.\qed\enddemo

In our situation $V=\P^8$, $\Cal L=\Cal O(dH)$ (where $d$ is the
degree of the curve $C$), and $\Cal E$ is the linear system
corresponding to the rational map $c=\alpha$.
We note that the support
$|E|$ of $E\hookrightarrow \P^8\times \P^N$ is described
set@-theoretically by
$$\multline
|E|=\{(\sigma,X)\in \P^8\times\P^N :\\
\text{$X$ is a limit of $\alpha(\sigma(t))$ for some curve germ
$\sigma(t)\subset \P^8$ centered at $\sigma\in S$}\}\quad,
\endmultline$$
so that it records the behavior of $\alpha$ as one approaches its base
locus $S$. Since $E$ is identified with the projective normal cone of
$S$ in $\P^8$, it is a scheme of pure dimension~7; invariably this
will turn out to be reducible and nonreduced. Often challenging is the
computation of the multiplicities $m_i$ of the various components $E_i$
of $E$; for our specific problem all this information can be found in
\cite{A-F5},
and it will be recalled in the next section.
In \S3
we will compute explicit expressions
$$E_i=(\epsilon_1 h^N H+\dots +\epsilon_8 h^{N-7} H^8)\cap [\P^8\times
\P^N]\quad, $$
yielding
$$L_i=\sum_{k\ge 0}\left(\sum_{j=0}^k \frac{(-d)^{k-j}\epsilon_{j+1}}
{j!(k-j)!}\right) \frac{H^{k+1}}{k+1}\quad.$$
According to Proposition~1.1, the a.p.p.~can be
computed by expanding
$$\exp(d H)\cdot\big(1-(m_1 L_1+\dots+m_k L_k)\big)\quad.$$
This will be our main tool in \S4 and \S5.

\example{Example~1.1} As an illustration, we describe the components of
$E$ for $C$ a smooth curve of degree $d\ge 2$, with only ordinary
flexes. Recall (\cite{A-F2}) that in this case the base locus $S$
consists of the set of rank@-1 matrices whose image is a point of~$C$.
We will see (\S2)
that $E$ consists of one component dominating $S$, and components
dominating the set of matrices whose image is an inflection point of $C$.

More precisely, the first component is supported on the
locus $G\subset\P^8\times\P^N$:{\eightpoint
$$\multline
G=\overline{
\{(\sigma,C_\sigma)\,|\,\text{$\im\sigma\in C$, and $C_\sigma$ is the
union $\ell\cup c$ of a $(d-2)$@-fold line $\ell$}}\\
\overline{
\text{supported on $\ker\sigma$ and a 
nonsingular conic $c$ tangent to $\ell$}\}}.
\endmultline$$}
\noindent Computing the class of this locus is a standard exercise in
the enumerative geometry of conics, and we obtain
$$[G]=6d\,H^5 h^{N-4}+4d (5d-9)\, H^6 h^{N-5}+6d(d-2)(5d-8)\,
H^7h^{N-6}\quad,$$
and the corresponding class
$$L_G=\frac{d H^5}{20}-\frac{d(5 d+18) H^6}{360}+\frac{d(9d+8)H^7} {420}
-\frac{d^2 H^8}{60}$$
in $\P^8$. The multiplicity of this component in the projective normal
cone turns out to be $2$ (Fact 2(ii) in \S2).

For each flex $p$ on $C$ we will also find a component of $E$
supported on $F\subset \P^8\times\P^N$:{\eightpoint
$$\multline
F=\overline{
\{(\sigma,C_\sigma)|\text{$\im\sigma=p$, and $C_\sigma$ is the union
of a $(d-3)$@-fold line $\ell$ supported on $\ker \sigma$}}\\ 
\overline{
\text{and a cuspidal cubic $c$ with cuspidal tangent $\ell$}\}}.
\endmultline$$}
\noindent Again the computation of the class of this locus in
$\P^8\times\P^N$ is not hard,
and yields 
$$L_F=\frac{H^6}{144}-\frac{H^7}{70}+\frac{197 H^8}{13440}$$
in $\P^8$; the multiplicity of $F$ in the projective normal cone will
be found to be~3 (Fact 4(ii) in \S2).
Since a smooth curve of degree $d\ge 2$ (and only ordinary flexes) has
$3 d(d-2)$ flexes, the predegree polynomial of such a curve is,
according to Proposition~1.1,
$$\exp(d H)\cdot(1-2\cdot L_G- 3d(d-2)\, 3\cdot L_F)$$
{\eightpoint
$$\multline
=1+d\,H+d^2\,\frac{H^2}2+d^3\,\frac{H^3}{3!}+d^4\,\frac{H^4}{4!}+(d^5-12d)\,
\frac{H^5}{5!}+(d^6-97 d^2+162 d)\,\frac{H^6}{6!}\\
+(d^7-427 d^3+1566 d^2-1488 d)\,\frac{H^7}{7!}
+(d^8-1372 d^4+7992 d^3-15879 d^2+10638 d)\,\frac{H^8}{8!}\quad.
\endmultline$$}
The coefficient of $\frac{H^8}{8!}$ reproduces the result
of the computation in \cite{A-F2} for $d\ge 3$.
Also note that, for $d=2$, this expression reduces to
$$1 + 2 H + \frac{4H^2}2 + \frac{8H^3}{3!} + \frac{16H^4}{4!} +
\frac{8H^5}{5!}\quad,$$
the predegree polynomial for a smooth conic, in agreement with
\cite{A-F3}, \S4.2.
We note in passing that the expression does {\it not\/} yield the
a.p.p.~of a line for $d=1$; this is not surprising, since a line is
not a curve with ordinary flexes.
\endexample


\head \S2. Limits of plane curves---summary of results\endhead

In this section we recall the results from \cite{A-F5} which we need
for the enumerative computations in this paper.

As we saw in \S1, we are interested in the structure of the
projective normal cone $E$ of the base scheme $S$ of the rational map
$$c:\P^8 \dashrightarrow \P^N$$
extending the action of $\PGL(3)$ on a given plane curve $C$ of
degree $d$.  
Now $S\subset \P^8$ consists of all matrices whose image is contained
in $C$; in particular, $S$ has exactly one component for each
component of $C$. More precisely, if no component of $C$ is a line,
then
$$|S|\cong \P^2\times |C|\subset \P^2\times\P^2 \subset \P^8\quad:$$
$S$ consists of rank@-1 matrices with arbitrary kernel, and image a
point of $C$. Every linear component $\ell$ of $C$ contributes a
5@-dimensional component to $S$, consisting of the $\P^5$ of
rank@-$\le 2$ matrices whose image is contained in $\ell$.

We have realized $E$ set@-theoretically as a subset of pure
dimension~$7$ of $\P^8\times \P^N$:
$$\multline
|E|=\{(\sigma,X)\in \P^8\times\P^N :\\
\text{$X$ is a limit of $c(\sigma(t))$ for some curve germ
$\sigma(t)\subset \P^8$ centered at $\sigma\in S$}\}\quad.
\endmultline$$
We are interested in a description of the components of this
locus, as well as the multiplicities with which they appear in $E$.
A given component may arise in several ways according to the procedure
described in this section; its multiplicity in $E$ will be
understood to be the sum of all multiplicities listed in each case. 

A first rough description of the components of $E$ can be given in terms
of the locus on $S$ they dominate:

\proclaim{Fact 1} There is one component of $E$ dominating each
component of $S$ (hence, one for each component of $C$), and components
dominating loci $\cong\P^2$:
$$\{\sigma\in\P^8\,|\,\text{$\sigma$ is a rank@-1 matrix with image
$p\in C$}\}$$
where $p$ is either a flex or a singular point of $C$.\endproclaim
We call the first kind of components `global', and the second kind
`local'. 

Components are usually best described as orbit closures of specific
elements $(\sigma,C_\sigma)$ of $\P^8\times \P^N$ under the induced
(right) action of $\PGL(3)$.
In each case $C_\sigma$ will be the limit obtained along a germ
centered at $\sigma$; thus it will be clear a priori that the given
locus is a component of $E$. The content of the results listed below
is that they provide an exhaustive list of all components of $E$ for a
given curve, and compute the multiplicity with which each component
appears.
Also, of course in each case $C_\sigma$ will be a curve with small
linear orbit; these curves have been studied in \cite{A-F3} and
\cite{A-F4}, and we use the terminology employed there.

Global components are easy to describe precisely:
\proclaim{Fact 2} (i) Let $\ell$ be a line appearing with multiplicity
$m$ in $C$, and let $\lambda$ be the $(d-m)$@-tuple of points cut out
on $\ell$ by the other components of $C$. Then the component of $E$
corresponding to $\ell$ is the orbit closure of{\eightpoint
$$\multline
(\sigma,C_\sigma),\text{where $\sigma$ is a rank@-2 matrix with
image $\ell$, and $C_\sigma$ is a fan consisting of a star}\\
\text{centered at $\ker\sigma$ and reproducing projectively the
tuple $\lambda$, and of a residual $m$@-fold line}\quad,
\endmultline$$}
with multiplicity $m$.

(ii) Let $C'$ be a non@-linear component appearing with multiplicity
$m$ in $C$. Then the component of $E$ corresponding to $C'$ is the
closure of the locus{\eightpoint
$$\multline
\{(\sigma,C_\sigma) \in\P^8\times\P^N\,|\,\text{$\sigma$ is a rank@-1
matrix with image a point of $C'$, and $C_\sigma$ consists of a}\\
\text{$(d-2m)$@-fold line supported on $\ker\sigma$, and of an
$m$@-fold smooth conic tangent to $\ker\sigma$}\}\quad,
\endmultline$$}
with multiplicity $2m$.
\endproclaim

We call components as in part (i)
{\it components of type~I,\/} and components as in part (ii)
{\it components of type~II.\/}

Local components of $E$ are substantially harder to describe, since the
germs of curves $\sigma(t)$ in $\P^8$ giving rise to such components
have to be carefully tailored to the local features of $C$. As shown
in \cite{A-F5},
only
two kinds of germs must be considered, requiring separate discussions:
one kind (1@-parameter subgroups, or 1@-PS for short) accounts for
limits with multiplicative stabilizer; the other will be responsible
for limits with additive stabilizer.

We start with the (simpler) case of 1@-PS limits. Again, we first
give a rough description of the situation.
\proclaim{Fact 3} Let $p$ be either a flex or a singular point of $C$.
For each line in the tangent cone to $C$ at $p$, there is a
corresponding Newton polygon.
The possible components of $E$ due to 1@-PS centered at $p$ are
indexed by sides of these Newton polygons;
further, an additional
component is present if the tangent cone is supported on at least three
distinct lines.\endproclaim
To be more precise, suppose $p$ has multiplicity $m$, and denote by $\lambda$
the tangent cone to $C$ at $p$ (hence $\lambda$ determines an
$m$@-tuple in the pencil of lines through $p$). 
\proclaim{Fact 4(i)} The component present exactly when $\lambda$ is
supported on three or more distinct lines is the orbit closure
of{\eightpoint
$$\multline
(\sigma,C_\sigma),\text{where $\sigma$ is a rank@-1 matrix
whose image is $p$, and $C_\sigma$ is a fan consisting of a}\\
\text{star projectively equivalent to $\lambda$, and of a residual
$(d-m)$@-fold line supported on $\ker\sigma$}\quad,
\endmultline$$}
with multiplicity $m A$, where $A$ is the number of automorphisms of
$\lambda$ as a tuple in the pencil of lines through $p$.
\endproclaim
\noindent
(The reason why this locus is not a component of $E$ if $\lambda$ is
supported on $\le 2$ lines is simply that it is not big enough to be
one: it is immediately checked that this locus has dimension~7 if and
only if $\lambda$ is supported on $\ge 3$ lines.) We call such
components {\it components of type~III.\/}

To determine the components corresponding to a line $\ell$ in the
tangent cone, choose coordinates $(x:y:z)$ in $\P^2$ so that
$p=(1:0:0)$ and $\ell$ is the line $z=0$; then 
consider the Newton polygon for the curve, that is, the boundary of the
convex hull of the union of the positive quadrants with origin at
the points $(j,k)$ for which the coefficient of $x^i y^j z^k$ in
the equation for $C$ is nonzero (see \cite{B-K}, p.~380).
Note that the part of the Newton polygon consisting of line segments with 
slope strictly between $-1$ and $0$ does not depend on the choice
of coordinates.
Consider the 1@-PS
$$\sigma(t)=\pmatrix
1 & 0 & 0\\
0 & t^b & 0\\
0 & 0 & t^c\endpmatrix\quad,$$
with $1\le b< c$ relatively prime integers, and $-b/c$ a slope of a
side of the Newton polygon for $C$.
\proclaim{Fact 4(ii)} For each line $\ell$ in the tangent cone of $C$,
and for each 1@-PS selected by the above procedure, there is a component
$E'$ of $E$ supported on the orbit closure of
{\eightpoint$$
(\sigma,C_\sigma),\text{where $C_\sigma$ is the limit 
as $t\to0$ of $C$ along
the selected 1@-PS $\sigma(t)$, and $\sigma=\sigma(0)$}\quad,
$$}
provided this locus has dimension~7.
If $x^{\overline q} y^r z^q \prod_{j=1}^S\left(y^c+\alpha_j
x^{c-b} z^b\right)$ is the limit obtained along the 1@-PS $\sigma(t)$,
then the contribution to the multiplicity of $E'$ is
$$(S b c+r b+q c)\,\frac A\delta\quad,$$
where $A$ is the number of components of the stabilizer of the limit,
and $\delta$ is the degree of the map from 
$E'$ to its image in $\P^N$.\endproclaim
The limits appearing in this statement are among the curves with
small orbit studied in \cite{A-F3}.
The number $\delta$ is 1 unless $c=2$ and $q=\overline q$, in which case
it is $2$ (see \cite{A-F5}).
The number $A/\delta$ can be computed directly in terms of
the tuple $\{\alpha_j\}$ (see \cite{A-F3}, Lemma 3.1).
We will see in \S3
that this factor is absorbed by
other terms in the computation of the contribution of such components.

We call components arising as in Fact~4(ii)
{\it components of type~IV.\/}

In order to visualize part of this somewhat complicated recipe, note
that if $(j_0,k_0)$, $(j_i,k_1)$, $j_0<j_1$, are vertices of a side of
the Newton polygon of $C$ of slope strictly between $-1$ and $0$,
then the corresponding multiplicity (provided the locus specified in the
statement has dimension 7) is
$$\frac{j_1 k_0-j_0 k_1}S \frac A\delta\quad,$$
where $S+1$ is the number of lattice points on the selected side.
Also, note that $\overline q=d-j_1-k_1$, $r=j_0$, and $q=k_1$ with
these notations; and $\delta=2$ exactly when $(j_0,k_0)$, $(j_i,k_1)$,
and $(d,0)$ lie on a line with slope $-1/2$.
The tuple $\{\alpha_j\}$ is determined by the specific coefficients
appearing along the side.

\example{Example~2.1} 
Suppose that $C$ has a general multiple point at $p$, by which we mean
an ordinary multiple point such that
the tangent line to each branch intersects that
branch with multiplicity $2$ at $p$. 
Let $m$ be the multiplicity of $C$ at $p$.
For each line in the tangent cone, the Newton
polygon contains exactly one side as in the prescription given above,
from $(m-1,1)$ to $(m+1,0)$; each line then contributes a multiplicity
of $(m+1)A/\delta$ to the component consisting of the orbit closure of
{\eightpoint$$
(\sigma,C_\sigma),\text{where $\tsize\sigma=\pmatrix 1 & 0 & 0\\
0 & 0 & 0\\ 0 & 0 & 0\endpmatrix$, and $C_\sigma$ is the curve
$x^{d-m-1} y^{m-1} (y^2+xz)=0$}\quad.
$$}
This component therefore appears in $E$ with multiplicity $m(m+1)A/\delta$.
Note that here $\delta=2$ exactly when the curve has degree $m+1$.
Also, if $m\ge3$ 
we find one component supported on the orbit closure of{\eightpoint
$$\multline
(\sigma,C_\sigma),\text{where $\sigma$ is a rank@-1 matrix whose
image is $p$, and $C_\sigma$ is a fan consisting of a}\\
\text{translate of the tangent cone at $p$, and of a residual
$(d-m)$@-fold line supported on $\ker\sigma$}\quad,
\endmultline$$}
with multiplicity $m$. 
\endexample

The real subtleties in the discussion occur in the next and last case,
dealing with limits with additive stabilizer. The components of $E$
detect an interaction between different (formal) branches of $C$
sharing a tangent at a singular point.
This phenomenon does not occur for e.g., ordinary multiple points.

Consider a line in the tangent cone to $C$ at $p$, and as above
choose coordinates so that $p=(1:0:0)$, and the line is $z=0$. 
Let $m$ be the multiplicity of $C$ at $p$.
It is well-known (cf.~\cite{B-K}) that there are $m$ formal
branches of $C$ at $p$, where nonreduced branches 
are counted according to their
multiplicity. For a general choice of $y$, these can be written
$$z=f(y)=\sum_i \gamma_{\lambda_i} y^{\lambda_i}\quad,$$ 
where $f(y)$ is a power series with fractional exponents 
$\lambda_i\in\Q$, $\lambda_0<\lambda_1<\dots$.

Let $B$ be the collection of 
all $m$ branches of the curve at $p$.
We then have a finite sequence of rational numbers $c>1$,
determined as those numbers $c$ for which at least two of the branches
tangent to $z=0$ agree modulo $y^c$, differ at $y^c$, and
satisfy $\lambda_0<c$.
Call $B_c$ the collection of those branches. 

Each $c$ determines a finite number of
truncations $\underline{f(y)}$: these are the truncations at 
$y^c$ (excluding $y^c$) of the branches in $B_c$.
These truncations determine germs
$$\sigma(t)=\pmatrix 1 & 0 & 0\\
t^a & t^{ab} & 0\\
\underline{f(t^a)} & \underline{f'(t^a) t^{ab}} & t^{ac}
\endpmatrix\quad,$$
where $b=\frac{c-\lambda_0}2+1$, and $a$ is the least positive integer
clearing all denominators in the exponents. 
We identify truncations if the corresponding germs are equivalent
after reparametrization, that is,
after multiplication on the right by $\pmatrix 1 & 0 & 0\\
0 & \eta^{ab} & 0\\
0 & 0 & \eta^{ac}\endpmatrix$, with $\eta$ a primitive $a$@-th root
of unity.

To each such germ we associate two numbers $\ell$ and $W$.
The number $\ell$ is defined as the least positive integer $\mu$
such that $\underline{f(y^{\mu})}$ has integer exponents.
The weight $W$ is defined as follows. For each branch
$\beta$ in $B$, let $v_{\beta}$ 
be the first exponent at which $\beta$ and $\underline{f(y)}$ differ,
and let $w_{\beta}$ be the minimum of $c$ and $v_{\beta}$. Then
$W$ is the sum $\sum w_{\beta}\,$.
\proclaim{Fact 5} Each germ $\sigma(t)$ contributes a component to
$E$: the orbit closure of {\eightpoint$$
(\sigma,C_\sigma),\text{where $C_\sigma$ is the limit of $C$ along
the germ $\sigma(t)$, and $\sigma=\sigma(0)$}\quad,
$$}
with multiplicity $\ell\, W A$, where $A$ is the number of components of
the stabilizer of~$C_\sigma$.
\endproclaim
The limits $C_\sigma$ appearing in this statement consist of unions of
quadritangent conics, plus possibly a multiple of the distinguished
tangent; these curves have been studied in \cite{A-F3}, \S4.1.
For enumerative purposes, they can be described in terms of the
multiplicities $s_i$ of the different conics, and of the number $A$ of
components of their stabilizer. As in the case of 1@-PS limits, this
number $A$ will be absorbed by other terms in the computation of the
contribution to the predegree of $C$.

We call the components identified in Fact~5
{\it components of type~V.}

An example will clarify the procedure described above.
\example{Example~2.2} Consider the quartic given in affine coordinates by
$$(y^2-xz)^2=y^3 z\quad.$$
Expanding at the origin gives two formal branches
$$z=y^2\pm y^{5/2}+\dots\quad;$$
with the notations used above: $c=\frac 52$, $b=\frac{5/2-2}2+1=\frac
54$, and $\underline{f(y)}=y^2$; hence the weight $W$ is $\frac
52+\frac 52=5$, $\ell=1$, and the germ determined by the truncation is
$$\sigma(t)=\pmatrix 1 & 0 & 0\\
t^4 & t^5 & 0\\
t^8 & 2 t^9 & t^{10}\endpmatrix\quad.$$
The corresponding component of $E$ is the orbit closure of
{\eightpoint$$
(\sigma,C_\sigma),\text{where $\tsize\sigma=\pmatrix 1 & 0 & 0\\
0 & 0 & 0\\ 0 & 0 & 0\endpmatrix$, and $C_\sigma$ is the curve
$({y^2} -x z + {x^2} )({y^2} -x z- {x^2})$}\quad;
$$}
one checks $A=4$,
and concludes that the multiplicity of this component in $E$ is
$1\cdot 5\cdot 4=20$.\endexample

To close the section, we remark that {\it not all\/} singular points
of (the support of) a curve contribute components to the projective
normal cone:
\example{Example~2.3} If $\ell_1$, $\ell_2$ are lines contained in $C$
(with any multiplicity), and $p=\ell_1\cap\ell_2$ is not a point of
the remainder of the curve, then $p$ does {\it not\/} contribute a
component to $E$.

Indeed, the tangent cone to $C$ at $p$ consists of only two lines, so
there are no components of type~III; next, the Newton polygon at $p$
with respect to either line has no sides of slope between $-1$ and
$0$, so there are no components of type~IV; finally, the branches of
$C$ at $p$ only consist of lines, so they do not interact in the sense
of providing a `truncation' as in Fact~5.
\endexample


\head \S3. Contributions to the adjusted predegree polynomial\endhead

The task in this section is to apply the results of \cite{A-F3},
\cite{A-F4}
and obtain explicit expressions for the contributions to the adjusted
predegree polynomials of a curve $C$ due to the various possible
components of the corresponding projective normal cone $E$. Together
with the description of the projective normal cone recalled in \S2,
the results of this section yield a procedure computing the predegree
polynomial of any given plane curve, in terms of the multiplicities of
its components and a description of its flexes and singular points.

Recall from \S1
that we have expressed the adjusted predegree polynomial (a.p.p.) of a
curve as
$$\exp(dH)\cdot(1-(m_1 L_1+\dots+m_k L_k))\quad;$$
our objective here is to obtain explicit expressions for the
different `correction' terms
$$-m_i L_i $$
due to the various components of the projective normal cone described
in \S2. The results will be used in \S4 and \S5
to obtain explicit expressions for contributions to the a.p.p.~due to
various features of a plane curve.
A correction term $-m_i L_i$ yields an {\it additive contribution\/}
$$\exp(d H)\cdot (-m_i L_i)$$
to the a.p.p.~of a curve of degree~$d$. All expressions $-m_i L_i$
will only have terms of degree $3$ or higher in $H$; those
corresponding to {\it local\/} components will only have terms of
degree $6$ or higher. Hence, the effect of a local correction term on
the a.p.p.~of a curve can also be expressed as a {\it multiplicative
contribution\/} by
$$(1-m_i L_i)\quad;$$
we will often prefer this alternative, since it does not involve the
degree of the curve. Also, sometimes we may list the effect of a
component as a correction term to the predegree of a curve, taking
account of other effects such as the number of flexes absorbed by a
given singularity.

In Propositions 3.1--3.5 below
we will compute the correction terms $-m_i L_i\,$.
As in \S2, we start with the global components.

\subheading{\S3.1. Type~I contributions}

\proclaim{Proposition~3.1} Let $\ell$ be a line appearing with
multiplicity $m$ in $C$, and let $r_i$ denote the multiplicities of
the intersections of $\ell$ with the rest of $C$. Then the correction
term due to $\ell$ is the antiderivative (w.r.t.~$H$) with 0 constant
term of
$$-\frac{m^3}2\,\exp(-d H)\,H^2\,\prod_i\left(1+r_i H+\frac{r_i^2 H^2}2
\right)\quad.$$
Explicitly:{\eightpoint
$$\multline
-\left(\frac{m^3 H^3}6-\frac{m^4 H^4}8 + \frac{m^5 H^5}{20} -
\frac{m^3(m^3+\sum r_i^3) H^6}{72}+\frac{m^3(m^4+4 m\sum r_i^3+3 \sum
r_i^4)H^7}{336}\right.\\
\left.- \frac{m^3(m^5+10 m^2 \sum r_i^3+15 m \sum r_i^4+6 \sum
r_i^5)H^8}{1920}\right)\quad.
\endmultline$$}
\endproclaim

\demo{Proof}
According to Fact~2(i) in \S2,
the component $E_\ell$ of $E$ corresponding to $\ell$ is the orbit
closure in $\P^8\times\P^N$ of $(\sigma,C_\sigma)$, where $\sigma$ has
image $\ell$ and $C_\sigma$ is a fan consisting of an $m$@-fold line
and a star $C'_\sigma$ of lines with multiplicities $r_1$, $r_2$, etc.,
centered at $\ker\sigma$. Denote by  
$$[E_\ell]=(\epsilon_1 H h^N+\dots +\epsilon_8 H^8 h^{N-7})\cap
[\P^8\times \P^N]$$
the class of this component, so that $\epsilon_i=H^{8-i} h^{i-1} \cdot
[E_\ell]$. 
\proclaim{Claim} Let $\beta_0+\beta_1 H+\dots+\beta_5 H^5$ be the
adjusted predegree polynomial of $C'_\sigma$. Then
$$\epsilon_i=\left\{\aligned
0\quad\qquad & i<3\\
\frac{m^2}2(i-1)!\,\beta_{i-3} \quad & i\ge 3
\endaligned\right.$$
\endproclaim
To see this, consider the embedding
$$\P^{N'}\times\P^2 @>>> \P^N\quad,$$
where $\P^{N'}$ parametrizes plane curves of degree $d-m$, $\P^2$
parametrizes lines, and the embedding attaches to a given curve of
degree $d-m$ an $m$@-fold line. We get an embedding
$$(\P^8\times\P^{N'})\times \P^2 @>\iota>> \P^8 \times \P^N\quad,$$
and it is readily understood that $E_\ell=\iota(E'_\ell \times \P^2),$
where $E'_\ell$ is the orbit closure of $(\sigma, C'_\sigma)$. Pulling
back to $(\P^8\times\P^{N'})\times \P^2$, we see that $\epsilon_i=0$
for $i<3$, and
$$\epsilon_i=m^2\binom{i-1}2 H^{8-i}{h'}^{i-3} \cdot [E'_\ell]$$
for $i\ge 3$, where $h'$ is the hyperplane in $\P^{N'}$.
Now note that $E'_\ell$ is the part of the closure of the graph of the
map
$$\P^8 \dashrightarrow \P^{N'}$$
(extending the action of $\PGL(3)$ on the star $C'_\sigma$) over the
$\P^5$ of matrices whose image is $\subset\ell$. By Remark~2.4 in
\cite{A-F4}
$$H^{8-i} {h'}^{i-3}\cdot [E'_\ell]=(i-3)!\,\beta_{i-3}\quad,$$
and the claim
follows.

The a.p.p.~for a star is computed in Theorem~2.5 in \cite{A-F4}:
$$\beta_0+\beta_1 H+\dots+\beta_5 H^5=\left\{\prod_i\left( 1+r_i H
+\frac{r_i^2 H^2}2\right)\right\}_5$$
($\{\}_5$ denotes truncation to $H^5$). Also, the multiplicity of this
component of $E$ is $m$, according to Fact~2(i) in \S2.
By the claim
and Proposition~1.1,
the correction term is therefore
$$-m\sum_{k\ge 0}\left(\sum_{j=0}^k \frac{(-d)^{k-j}\epsilon_{j+1}}
{j!(k-j)!}\right) \frac{H^{k+1}}{k+1} =-\frac{m^3}2 \sum_{k\ge 0}
\left(\sum_{j=2}^k \frac{(-d)^{k-j}}{(k-j)!}\beta_{j-2}\right)
\frac{H^{k+1}}{k+1}\quad,$$
yielding the expressions given in the statement.\qed\enddemo

\example{Example~3.1} The a.p.p.~of a curve consisting of a union of
lines, with multiplicity $m_i$ and no three meeting at a point, is
$$\prod_{i}\left(1+m_i H+\frac{m_i^2 H^2}2\right)$$
(by our notational convention, this expression stands for its
truncation at $H^8$).

Indeed, by Example~2.3
there are no components of $E$ due to the points of intersection of
such a configuration of lines; the only components are therefore those
corresponding to the lines themselves. Using Proposition~3.1,
the total correction term evaluates to
{\eightpoint
$$\multline
-\left(\frac{\sum m_i^3 H^3}6-\frac{\sum m_i^4 H^4}8 + \frac{\sum
m_i^5 H^5}{20} - \frac{(\sum m_i^3)^2 H^6}{72}+\frac{(7(\sum
m_i^3)(\sum m_i^4)-6\sum m_i^7)H^7}{336}\right.\\
\left.- \frac{(15(\sum m_i^4)^2+16 (\sum m_i^3)(\sum m_i^5)-30 \sum
m_i^8)H^8}{1920}\right)\quad.
\endmultline$$}
Applying Proposition~1.1
yields the expression given in the statement.

This computation reproduces results from \S2 of \cite{A-F4},
where a more general `multiplicativity' of adjusted predegree
polynomials for configurations of lines meeting transversally is
discussed.\endexample

\subheading{\S3.2. Type~II contributions}
Next, we consider nonlinear components of $C$:
\proclaim{Proposition~3.2} Let $C'$ be a component of $C$ of degree $e>1$,
appearing with multiplicity $m$ in $C$. Then the correction term due
to $C'$ is
$$-2e m^5\left(\frac{H^5}{20}-\frac{(5d+18 m) H^6}{360} + \frac{(9d+8
m)mH^7}{420}-\frac{dm^2 H^8}{60}\right)\quad.$$
\endproclaim
\demo{Proof} According to Fact 2(ii), the corresponding component of
$E$ is the locus $E_{C'}$ of $(\sigma,C_\sigma)$, where the image of
$\sigma$ is a point of $C'$ and $C_\sigma$ consists of a
$(d-2m)$@-fold line supported on $\ker\sigma$, and of an
$m$-fold conic tangent to $\ker\sigma$. Let
$$[E_{C'}]=(\epsilon_1 Hh^N +\dots +\epsilon_8  H^8h^{N-7})\cap
[\P^8\times \P^N];$$
then $\epsilon_i= H^{8-i}h^{i-1}\cdot [E_{C'}]$. To evaluate this,
note that $E_{C'}$ is contained in $B\times\P^N\subset \P^8\times\P^N$,
where $B=\P^2\times C'$ is the set of rank@-1 matrices $\sigma$ with
image on $C'$. Denoting by $k$ the pull@-back to $B$ of the hyperplane
class from the $\P^2$ factor, and by $\ell$ the pull@-back of the
restriction of the hyperplane class from the other factor, we have
$$\epsilon_i=(k+\ell)^{8-i}h^{i-1}\cdot [E_{C'}]=(8-i) k^{7-i}\ell
h^{i-1}\cdot [E_{C'}]\quad;$$
in particular $\epsilon_i=0$ unless $i=5$, $6$, or $7$. The class $\ell$
splits $E_{C'}$ into $e$ components, each of which consists
of points $(\sigma,C_\sigma)$ with $\sigma$ constrained to have a
fixed image. Also note that intersecting by $k$ amounts to imposing a
linear condition on the distinguished tangent line in $C_\sigma$;
therefore, $\epsilon_i=$ $(8-i)e$ times the number (counted with
multiplicity) of curves $C_\sigma$ through $i-1$
general points, with tangent line constrained to contain $7-i$ general
points, where $i=5$, $6$, or $7$.

For these values of $i$, 
the corresponding number of configurations (in case $d>2m$) is
computed by arguing as in \cite{A-F3}, Proposition~4.1:
$$\epsilon_i=(8-i)\, e\,\frac{(i-1)!}{6!}
\frac{\partial^{7-i}}{\partial \overline q^{7-i}} P(\overline
q)|_{\overline q=d-2m}\quad,$$
where $P(\overline q)$ is the polynomial giving the degree 
for a curve
such as $C_\sigma$, with distinguished tangent taken with multiplicity
$\overline q$. This is the coefficient of $t^6/6!$ in the a.p.p.~for
$C_\sigma$ (computed in \S4.2 of \cite{A-F3}:
set $n=2$; $m=\overline m=1$; $S=s_1=m$; $r=q=0$ in the formulas given there),
divided by $4$, the degree of the stabilizer:
$$P(\overline q)=12 m^5 \overline q+30 m^4\overline q^2\quad.$$
The same formula holds in the case $d=2m$. This yields
$$[E_{C'}]=em^4 \left(6H^5 h^{N-4}+ 4(5d-9m) H^6 h^{N-3}+
6(5d-8m)(d-2m) H^7 h^{N-2}\right)\quad.$$
According to Fact~2(ii) in \S2
this locus appears in $E$ with multiplicity $2m$. From this we obtain
the stated correction term.\qed\enddemo

\example{Example~3.2} If $C$ is reduced and irreducible, then the only
component of type~II considered in Proposition~3.2
is the one dominating the whole curve. Setting $e=d$, $m=1$ we get a
correction term of
$$-2d\left(\frac{H^5}{20}-\frac{(5 d+18) H^6}{360}+\frac{(9d+8)H^7}
{420} -\frac{d H^8}{60}\right)$$
agreeing with the class $-2\,L_G$ used in Example~1.1
\endexample

\subheading{\S3.3. Type~III contributions}
Moving on to the correction terms due to local features of the curve,
we first establish a technical lemma, which will be used in the proofs of
the statements that follow, and which explains a recurrent feature of
the correction terms we will compute.

The components of type~III, IV, and V, arising from
local features of the curve, consist of orbit closures of points
$(\sigma,C_\sigma)\in \P^8\times \P^N$, where $\sigma$ is a rank@-one
matrix with a given image point and $C_\sigma$ is a curve with a
distinguished line, that is supported on $\ker\sigma$ and has
multiplicity $\overline q=d-\rho$ (where $\rho$ changes from case to
case). Let $P(\overline q)$ denote the coefficient of $H^7$ in the
predegree polynomial for such a curve; this is always a polynomial of
degree at most two in $\overline q$. 
Also, let $\delta$ be the degree of the map from the component to its
image in $\P^N$. As pointed out already in \S2,
this number is $1$ in almost all cases.
\proclaim{Lemma~3.3.1} The corresponding contribution to the correction
term is
$$-\delta\left(\frac{P''(-\rho) H^6}{42\cdot 6!}+\frac{P'(-\rho) H^7}
{7\cdot 7!}+\frac{P(-\rho) H^8}{8!}\right)\quad.$$
\endproclaim
\demo{Proof}Let $E'$ denote a component of $E$ arising from a point
$p$ of the curve, and let 
$$[E']= (\epsilon_1 Hh^N +\dots +\epsilon_8  H^8h^{N-7})\cap
[\P^8\times \P^N]$$
be its class. Since $E'$ is the orbit closure of a point
$(\sigma,C_\sigma)\in\P^8\times\P^N$, with $\sigma$ a rank@-1 matrix
with image $p$, $E'$ is in fact contained in $\P^2\times\P^N \subset
\P^8\times \P^N$, where $\P^2$ consists of all rank@-1 matrices with
image $p$. If $k$ denotes the hyperplane class in $\P^2$, pulling back
to $\P^2\times\P^N$ shows that
$$\epsilon_i=k^{8-i} h^{i-1}\cdot [E']\quad;$$
this gives immediately $\epsilon_i=0$ unless $i=6$, $7$, or $8$. Also,
note that under the identification of $\P^2$ with rank@-1 matrices
$\sigma$ with fixed image, the class $k$ imposes a linear condition on
the line $\ker\sigma$. Now, $C_\sigma$ consists in each case of a
curve with a distinguished line supported on $\ker\sigma$, appearing
with multiplicity $\overline q=d-\rho$ in our notations. Let $P(\overline
q)=\alpha \overline q^2+\beta \overline q+\gamma$ be the polynomial in
$\overline q$ giving the coefficient of $H^7$ in the predegree
polynomial for such a curve.
Using Proposition~4.1 in
\cite{A-F3}
we get
$$\frac{\epsilon_i}{\delta}=\left\{\aligned
\tfrac{P''(d-\rho)}{42}\quad &i=6\\
\tfrac{P'(d-\rho)}{7}\quad &i=7\\
P(d-\rho)\quad &i=8\endaligned\right.\quad\text{and therefore}$$
$$[E']=\frac{2\alpha}{42}H^6
h^{N-5}+\frac{2\alpha(d-\rho)+\beta}{7} H^7 h^{N-6} 
+ (\alpha(d-\rho)^2+\beta(d-\rho)+\gamma) H^8 h^{N-7}.$$
Computing the corresponding correction term as prescribed in \S1
gives the stated expression.\qed\enddemo
This observation explains why the degree $d$ of $C$ does {\it not\/}
appear explicitly in the correction terms we will list.
Note that a similar phenomenon
also occurs in the second formula in Proposition~3.1.

Let $p$ be a singular point of $C$. As recalled in \S2, Fact~4(i),
a component of type~III of the projective normal cone is present if
the tangent cone to $C$ at $p$ is supported on $\ge 3$ distinct lines.
\proclaim{Proposition~3.3} Let $e_i$ denote the elementary symmetric
functions in the multiplicities of the distinct lines in the tangent
cone to $C$ at $p$ (so $e_1=$ the multiplicity of $C$ at $p$). Then
the correction term corresponding to this component is
$$-e_1(e_2 e_3-e_1 e_4-e_5)\left(\frac{H^6}{24}-\frac{e_1 H^7}{28}+
\frac{e_1^2 H^8}{64}\right)\quad.$$
\endproclaim
Note that the expression given in this statement vanishes
automatically if the tangent cone is supported on $\le 2$ lines.
\demo{Proof} Using Fact~4(i)
and Lemma~3.3.1,
the main ingredient in the computation is the polynomial $P(\overline
q)$ expressing the degree for a fan $C_\sigma$ with star projectively
equivalent to the tangent cone to $C$ at $p$, and residual
$\overline q$@-fold line. From \cite{A-F4}, Theorem~2.5(ii), this
polynomial is
$$P(\overline q)=\frac{630\, \overline q^2}{A}\left(e_2e_3-e_1
e_4-e_5\right)\quad,$$ 
where $A$ is the number of automorphisms of the tuple determined by
the lines in the tangent cone as elements of the pencil of lines
through $p$. By Lemma~3.3.1,
with $\overline q=d-e_1$, the correction term is
$$-\frac{(e_2 e_3-e_1 e_4-e_5)}{A}\left(\frac{H^6}{24}
-\frac{e_1 H^7}{28}+\frac{e_1^2 H^8}{64}\right)$$
times the multiplicity with which the component appears in the
projective normal cone. By Fact~4(i) this multiplicity is $e_1 A$,
and the statement follows.\qed\enddemo

\example{Example~3.3} If the tangent cone consists of $m$ distinct reduced
lines, then Proposition~3.3
evaluates its corresponding correction term as
$$-m\left(\binom m2 \binom m3-m \binom m4-\binom m5\right)
\left(\frac{H^6}{24}-\frac{m H^7}{28}+ \frac{m^2 H^8}{64}\right)
\quad,$$
that is
$$-m^2 (m-1)(m-2)(m^2+3 m-3)\left(\frac{H^6}{720}-\frac{m H^7}{840}+
\frac{m^2 H^8}{1920}\right) \quad.$$
As an illustration, consider a star of $d$ reduced lines through a
point. The point will contribute as above, with $m=d$; also, according
to Proposition~3.1
each line contributes{\eightpoint
$$\multline
-\left(\frac{H^3}6-\frac{H^4}8 + \frac{H^5}{20} -
\frac{(1+(d-1)^3) H^6}{72}+\frac{(1+4 (d-1)^3+3 (d-1)^4)H^7}{336}\right.\\
\left.- \frac{(1+10 (d-1)^3+15 (d-1)^4+6 (d-1)^5)H^8}{1920}\right)\quad.
\endmultline$$}
{}From the discussion of \S2,
we know that there are no other correction terms. Putting 
everything together and using Proposition~1.1,
the a.p.p.~of this curve is{\eightpoint
$$\multline
\exp(d H)\left(1-\frac{d H^3}{6}+\frac{d H^4}{8}-\frac{d H^5}{20}
-\frac{d^2(d-3)(d^3+3 d^2-11 d+12)H^6}{720}\right.\\
\left.+\frac{d^3(2 d^4-35 d^2+70 d-42)H^7}{1680} -\frac{d^4(d^4-16 d^2+30
d-16)H^8}{1920}\right)
\endmultline$$}
that is{\eightpoint
$$1+d H+\frac{d^2
H^2}2+\frac{d(d-1)(d+1)H^3}6+\frac{d(d-1)(d^2+d-3)H^4} {24}
+\frac{d(d-1)(d-2)(d^2+3d-3)H^5}{120}$$}
Note that the polynomial detects that the orbit closure of this curve
has dimension $\le 5$; of course the stated expression is the truncation 
$$\left\{\left(1+H+\frac{H^2}{2}\right)^d\right\}_5$$
as prescribed by \cite{A-F4}, Theorem~2.5(i).
In fact, Propositions~3.1 and 3.3
suffice to compute the a.p.p.~for an arbitrary configuration of lines
in the plane, recovering Theorem~2.8 in \cite{A-F4}.
\endexample

\subheading{\S3.4. Type~IV contributions}
Next, let $p$ be a singular or an inflection point of (the support of)
$C$,
and consider a line $\ell$ of the tangent cone to $C$ at $p$. We have
recalled in Fact~4(ii)
that these choices determine a Newton polygon, and that there are
components (of type~IV) of the projective normal cone corresponding to
the sides of this polygon of slope strictly between $-1$ and $0$.

Consider then such a side $\Sigma$, from $(j_0,k_0)$ to $(j_1,k_1)$
($j_0<j_1$); and let $S+1$ be the number of lattice points on $\Sigma$.
Let $\gamma_0$, \dots, $\gamma_S$ be the coefficients on $\Sigma$ of the
equation for $C$, and consider the $S$@-tuple in $\P^1$
determined by the polynomial
$$\gamma_0 \xi^S +\gamma_1 \xi^{S-1}\eta +\dots+\gamma_S
\eta^S\quad;$$
let $s_i$ be the multiplicities of the points of this $S$@-tuple (so
for example $S=\sum s_i$).

The side $\Sigma$ determines then the following expressions:

$\bullet$ $R(\Sigma)=(j_1 k_0-j_0 k_1)$, that is, twice the area of
the triangle with vertices at $(0,0)$, $(j_0,k_0)$, and $(j_1,k_1)$;

$\bullet$ a polynomial $G(\Sigma)=\frac 1S\left(4\sum_i
s_i^5\frac{H^6}{6!}-36 \sum_i s_i^6 \frac{H^7}{7!}+192\sum_i
s_i^7\frac{H^8}{8!}\right)$;

$\bullet$ and a polynomial $L(\Sigma)$ given by
{\eightpoint $$
\multline
\left(6 j_0^2 k_0^2 + 3 j_0 j_1 k_0^2 + j_1^2 k_0^2 + 3 j_0^2 k_0 k_1 + 
  4 j_0 j_1 k_0 k_1 + 3 j_1^2 k_0 k_1 + j_0^2 k_1^2 + 3 j_0 j_1 k_1^2 + 
  6 j_1^2 k_1^2\right)\frac{H^6}{6!} \\
-\left(30 j_0^3 k_0^2 + 18 j_0^2 j_1 k_0^2 + 9 j_0 j_1^2 k_0^2 + 
  3 j_1^3 k_0^2 + 30 j_0^2 k_0^3 + 12 j_0 j_1 k_0^3 + 
  3 j_1^2 k_0^3 + 12 j_0^3 k_0 k_1 + 18 j_0^2 j_1 k_0 k_1 \right. \\
+ 18 j_0 j_1^2 k_0 k_1 + 12 j_1^3 k_0 k_1 + 18 j_0^2 k_0^2 k_1
+ 18 j_0 j_1 k_0^2 k_1 + 9 j_1^2 k_0^2 k_1 + 3 j_0^3 k_1^2 + 
  9 j_0^2 j_1 k_1^2 + 18 j_0 j_1^2 k_1^2\\
\left.+ 30 j_1^3 k_1^2 + 9 j_0^2 k_0 k_1^2 + 18 j_0 j_1 k_0 k_1^2 + 18
j_1^2 k_0 k_1^2 + 3 j_0^2 k_1^3 + 12 j_0 j_1 k_1^3 + 30 j_1^2
k_1^3\right)\frac{H^7}{7!} \\
+\left(90 j_0^4 k_0^2 + 60 j_0^3 j_1 k_0^2 + 36 j_0^2 j_1^2 k_0^2 + 
  18 j_0 j_1^3 k_0^2 + 6 j_1^4 k_0^2 + 180 j_0^3 k_0^3 + 
  90 j_0^2 j_1 k_0^3 + 36 j_0 j_1^2 k_0^3\right.\\
+ 9 j_1^3 k_0^3 + 
  90 j_0^2 k_0^4 + 30 j_0 j_1 k_0^4 + 6 j_1^2 k_0^4 + 
  30 j_0^4 k_0 k_1 + 48 j_0^3 j_1 k_0 k_1 + 54 j_0^2 j_1^2 k_0 k_1 + 
  48 j_0 j_1^3 k_0 k_1\\
+ 30 j_1^4 k_0 k_1 + 90 j_0^3 k_0^2 k_1 + 
  108 j_0^2 j_1 k_0^2 k_1 + 81 j_0 j_1^2 k_0^2 k_1 + 
  36 j_1^3 k_0^2 k_1 + 60 j_0^2 k_0^3 k_1 + 48 j_0 j_1 k_0^3 k_1\\
+ 18 j_1^2 k_0^3 k_1 + 6 j_0^4 k_1^2 + 18 j_0^3 j_1 k_1^2 + 
  36 j_0^2 j_1^2 k_1^2 + 60 j_0 j_1^3 k_1^2 + 
  90 j_1^4 k_1^2 + 36 j_0^3 k_0 k_1^2 + 81 j_0^2 j_1 k_0 k_1^2\\
+ 108 j_0 j_1^2 k_0 k_1^2 + 90 j_1^3 k_0 k_1^2 + 
  36 j_0^2 k_0^2 k_1^2 + 54 j_0 j_1 k_0^2 k_1^2 + 
  36 j_1^2 k_0^2 k_1^2 + 9 j_0^3 k_1^3 + 36 j_0^2 j_1 k_1^3 + 
  90 j_0 j_1^2 k_1^3\\
\left. + 180 j_1^3 k_1^3 + 18 j_0^2 k_0 k_1^3 + 
  48 j_0 j_1 k_0 k_1^3 + 60 j_1^2 k_0 k_1^3 + 6 j_0^2 k_1^4 + 
  30 j_0 j_1 k_1^4 + 90 j_1^2 k_1^4\right)\frac{H^8}{8!}
\endmultline
$$}
This polynomial is symmetric in the vertices of $\Sigma$;
unfortunately, we do not have a more intrinsic interpretation for it.

\proclaim{Proposition~3.4} The correction term due to the selected line
$\ell$ in the tangent cone to $C$ at $p$ is
$$-\sum_{\Sigma} R(\Sigma)\left(L(\Sigma)-G(\Sigma)\right)\quad.$$
\endproclaim

\demo{Proof} This follows from Lemma~3.3.1
and Fact~4(ii).
Using the notations of Fact~4(ii),
for each side $\Sigma$ we need the coefficient of the term of degree
$7$ in the predegree polynomial for limit curves $C_\sigma$ with
equation 
$$x^{\overline q} y^r z^q \prod_{j=1}^S\left(y^c+\alpha_j x^{c-b}
z^b\right)\quad,$$
where 
$$\gamma_0 \xi^S +\gamma_1 \xi^{S-1}\eta +\dots+\gamma_S \eta^S=
\gamma_0\prod_j(\xi-\alpha_j \eta)\quad.$$
These are precisely the curves studied in \cite{A-F3};
the predegree polynomial for such curves is computed in Theorem~1.1 of
\cite{A-F3}.
In our situation, we have
$$r=j_0,\quad q=k_1,\quad \overline q=d-(j_1+k_1)$$
(hence we use $\rho=j_1+k_1$ when applying Lemma~3.3.1),
and
$$b=\frac{k_0-k_1}S\quad,\quad c=\frac{j_1-j_0}S\quad;$$
applying Lemma~3.3.1
to the polynomial in $\overline q$ obtained from Theorem~1.1 in \cite{A-F3}
gives the expression
$$-\frac {S\delta}A\left(L(\Sigma)-G(\Sigma)\right)\quad,$$
where $A$ denotes the number of components of the stabilizer of
$C_\sigma$,
and $\delta$ is as in Lemma~3.3.1.

According to Fact~4(ii),
the contribution to the multiplicity of this component due to $\Sigma$
is
$$(S b c+r b+q c)\frac A\delta=\frac{j_1 k_0-j_0 k_1}S \frac A\delta=
R(\Sigma)\,\frac A{S\delta}\quad;$$
the correction term is therefore as stated.\qed\enddemo

\example{Example~3.4} Suppose $p$ is a $k$@-flex of $C$, that is, a
nonsingular point of $C$ at which $C$ and its tangent line $\ell$ meet
with multiplicity $k$. (For example, an ordinary inflection point of
$C$ is a $3$@-flex in this terminology). The Newton polygon at $\ell$
has only one side $\Sigma$ with slope between $-1$ and $0$, with
vertices $(0,1)$ and $(k,0)$. We have $S=1$, and the expressions given
above evaluate to
$$\gather
R(\Sigma)=k\quad,\quad G(\Sigma)=\frac{4H^6}{6!}-\frac{36H^7}{7!}
+\frac{192H^8}{8!}\\
L(\Sigma)=\frac{k^2H^6}{6!}-\frac{(3k^2+3k^3)H^7}{7!}
+\frac{(6k^2+9k^3+6k^4)H^8}{8!}\quad,
\endgather$$
giving a correction term of
$$k(k-2)\left(\frac{(k+2) H^6}{720}-\frac{(k^2+3k+6)H^7}{1680}+
\frac{(2k^3+7k^2+16k+32)H^8}{13440}\right)\quad.$$
For $k=3$, this recovers the term $L_F$ used in Example~1.1.
\endexample

The analysis presented up to this point suffices already to compute
the predegree of an arbitrary plane curve with ordinary multiple
points; this case is analyzed in~\S4.

\subheading{\S3.5. Type~V contributions}
We are left with the case of components of the projective normal cone
$E$ of type~V, arising from the interaction of different formal
branches with the same tangent line at a point $p$ of $C$. As pointed
out in \S2,
contributions corresponding
to these components arise from truncations of power
series with fractional exponents representing the different branches:
roughly, a contribution arises when two branches agree up to a certain
exponent $c$, but differ at that exponent. Truncating there determines
a germ $\sigma(t)$, centered at $\sigma=\sigma(0)$, and a limit
$C_\sigma$; the corresponding component consists of the orbit closure
of $(\sigma, C_\sigma)$. Further, the germ determines two
numbers $\ell$, $W$ (see Fact~5 in \S2).

Limits $C_\sigma$ obtained by this procedure consist of unions of
4@-tangent conics, and a multiple of the distinguished tangent, which
is supported on $\ker\sigma$. We let $s_i$ denote the multiplicities
with which the conics appear in $C_\sigma$, and write $S=\sum s_i$.
\proclaim{Proposition~3.5} With notations as above, the corresponding
correction term is
$$-\ell W\left(\frac{4(S^5-\sum_i s_i^5) H^6}{6!} -\frac{36(S^6-\sum_i
s_i^6) H^7}{7!}+\frac{192(S^7-\sum_i s_i^7)H^8}{8!}\right)\quad.$$
\endproclaim
\demo{Proof} This is obtained from Lemma~3.3.1
and Fact~5 in \S2,
by the procedure applied in Propositions~3.3 and 3.4.
The main ingredient is the predegree of the curves $C_\sigma$, which
is given in \S4.1 of \cite{A-F3}.
\qed\enddemo

\example{Example~3.5} As an illustration, we take the origin $(1:0:0)$ in
the curve
$$(y^2-x z)^2=y^3 z\quad.$$
As seen in Example~2.2,
only one truncation needs to be considered for this point; the
corresponding limit is a pair of distinct conics; further, $\ell=1$
and $W=5$. With notations as above we have $s_1=s_2=1$, so according
to Proposition~3.5
the corresponding correction term is
$$-5\left(\frac{H^6}6-\frac{31 H^7}{70}+\frac{3 H^8}5\right)\quad.$$
Applying Proposition~1.1,
this yields a contribution to the a.p.p.~of
$$-\left(\frac{5H^6}6+\frac{47 H^7}{42}+\frac{17 H^8}{21}\right)\quad;$$
in particular, the contribution due to this limit to the predegree of
the curve is
$$-8!\,\frac{17}{21}=-5\cdot 6528\quad.$$
This example belongs to a class of singular points which can be realized
on a quartic curve, and are analytically isomorphic to the singularity
$z^2=y^k$, $k=5$ (as in this example), $6$, $7$, or $8$. The
corresponding contribution to the predegree of the quartic turns out
to be $-k\cdot 6528$ in all cases (cf.~Example~5.4).
\endexample

\remark{Remark} As an immediate application of the results obtained
above, we can measure the effect on the contribution of a point $p$
due to taking a `multiple' of the curve on which $p$ lies.

If $C$ has ideal $(F(x:y:z))$ and $m$ is a positive integer, we let
$mC$ denote the curve with ideal $(F^m)$. Let $p\in C$, and assume the
contribution of $p$ to the a.p.p.~of $C$ is $K(H)$.
\proclaim{Claim} Then the contribution of $p$ to $mC$ is $K(m
H)$. \endproclaim
\demo{Proof} This follows from the homogeneity of the various
correction terms. The effect of replacing $C$ by $mC$ is that of
replacing $e_i$ by $m^i e_i$ in correction terms of type~III;
$(j_i,k_i)$ by $(m j_i,m k_i)$, $W$ by $mW$, and $S$, $\sum s_i^5$,
$\sum s_i^6$, $\sum s_i^7$ by $mS$, $m^5\sum s_i^5$, $m^6 \sum s_i^6$,
$m^7 \sum s_i^7$ respectively in correction terms of type~IV and
V. The claim
follows.\qed\enddemo

A similar homogeneity holds for global correction terms as well, so
that if $P(H)$ is the a.p.p.~of a curve $C$, then $P(m H)$ is the
a.p.p.~of its multiple $mC$. This can also be deduced by considering
the map $\P^{d(d+3)/2} @>>> \P^{md(md+3)/2}$ defined by $C \mapsto m
C$, a projection of the $m$@-th Veronese embedding.\endremark

\subheading{\S3.6. Summary} The results obtained in this section,
together with the discussion in \S2,
give an algorithm to compute the adjusted predegree polynomial of an
arbitrary plane curve. This will be illustrated in \S4 and \S5
by applying it to several classes of curves.

For reference we list here the contributions to the predegree of a
curve (with orbit of dimension~8) due to its features. Each of these
is obtained by applying Proposition~1.1
to the results obtained in Propositions~3.1--3.5,
obtaining corresponding additive contributions to the a.p.p., then
reading the coefficient of $\frac{H^8}{8!}$.

Assume $C$ has degree $d$. The predegree of its orbit closure is
obtained then by subtracting various contributions from $d^8$, indexed
here according to the corresponding type:

(I) A line appearing in $C$ with multiplicity $m$, meeting the rest of
the curve along a $(d-m)$@-tuple of points with multiplicities $r_i$,
gives a contribution of{\eightpoint
$$\multline {m^3} \left( {d^3} \left( 10 {d^2} - 15 d m +
6 {m^2} \right)  + 10
\left( 28 {d^2} - 48 d m + 21 {m^2} \right) \left({{\left( d - m
\right) }^3}- {\sum r_i^3} \right)\right.\\
\left.  -  45 \left( 8 d - 7 m \right) \left( {{\left( d - m \right) }^4}
- {\sum r_i^4} \right)  + 126 \left( {{\left( d - m \right) }^5} - {\sum
r_i^5} \right)  \right)\quad.
\endmultline$$}

(II) A component of $C$ of degree $e>1$ and appearing with
multiplicity $m$
contributes
$$16 d e m^5\left(7 d^2-18 d m+12 m^2\right)\quad.$$

Points $p\in C$ may contribute different terms:

(III) Let $e_i$ be the elementary symmetric functions in the
multiplicities of the distinct lines in the tangent cone to $C$ at
$p$. Then the corresponding contribution is
$$30 e_1(e_2 e_3-e_1 e_4-e_5)\left(28 d^2-48 d e_1+21 e_1^2\right)\quad.$$
(In particular, no such contribution is present if the tangent cone
consist of $<3$ distinct lines.)

(IV)
Let $\ell$ be a line of the tangent cone of $C$ at $p$, and let
$\Sigma$ denote the sides of slope strictly between $-1$ and $0$ of
the corresponding Newton polygon. With notations as in Proposition~3.4,
the contribution due to each $\Sigma$ is obtained by adding
$$\frac{16(j_1-k_0)}S \left(7 d^2 \sum s_i^5-18 d \sum s_i^6+12 \sum
s_i^7\right)$$
and{\eightpoint
$$\multline
(j_1 k_0 - j_0 k_1 )
  \left( 90 j_0^4 k_0^2 + 180 j_0^3 k_0^3 + 90 j_0^2 k_0^4 + 
    60 j_0^3 k_0^2 j_1 + 90 j_0^2 k_0^3 j_1 + 30 j_0 k_0^4 j_1 + 
    36 j_0^2 k_0^2 j_1^2\right.\\
+ 36 j_0 k_0^3 j_1^2 + 
    6 k_0^4 j_1^2 + 18 j_0 k_0^2 j_1^3 + 9 k_0^3 j_1^3 + 
    6 k_0^2 j_1^4 - 240 j_0^3 k_0^2 d - 240 j_0^2 k_0^3 d - 
    144 j_0^2 k_0^2 j_1 d\\
- 96 j_0 k_0^3 j_1 d - 
    72 j_0 k_0^2 j_1^2 d - 24 k_0^3 j_1^2 d - 
    24 k_0^2 j_1^3 d + 168 j_0^2 k_0^2 d^2 + 
    84 j_0 k_0^2 j_1 d^2 + 28 k_0^2 j_1^2 d^2 + 
    30 j_0^4 k_0 k_1\\
+ 90 j_0^3 k_0^2 k_1 + 60 j_0^2 k_0^3 k_1 + 
    48 j_0^3 k_0 j_1 k_1 + 108 j_0^2 k_0^2 j_1 k_1 + 
    48 j_0 k_0^3 j_1 k_1 + 54 j_0^2 k_0 j_1^2 k_1 + 
    81 j_0 k_0^2 j_1^2 k_1\\
+ 18 k_0^3 j_1^2 k_1 + 
    48 j_0 k_0 j_1^3 k_1 + 36 k_0^2 j_1^3 k_1 + 30 k_0 j_1^4 k_1 - 
    96 j_0^3 k_0 d k_1 - 144 j_0^2 k_0^2 d k_1 - 
    144 j_0^2 k_0 j_1 d k_1\\
- 144 j_0 k_0^2 j_1 d k_1 - 
    144 j_0 k_0 j_1^2 d k_1 - 72 k_0^2 j_1^2 d k_1 - 
    96 k_0 j_1^3 d k_1 + 84 j_0^2 k_0 d^2 k_1 + 
    112 j_0 k_0 j_1 d^2 k_1\\
+ 84 k_0 j_1^2 d^2 k_1 + 
    6 j_0^4 k_1^2 + 36 j_0^3 k_0 k_1^2 + 
    36 j_0^2 k_0^2 k_1^2 + 18 j_0^3 j_1 k_1^2 + 
    81 j_0^2 k_0 j_1 k_1^2 + 54 j_0 k_0^2 j_1 k_1^2 + 
    36 j_0^2 j_1^2 k_1^2\\
+ 108 j_0 k_0 j_1^2 k_1^2 + 
    36 k_0^2 j_1^2 k_1^2 + 60 j_0 j_1^3 k_1^2 + 
    90 k_0 j_1^3 k_1^2 + 90 j_1^4 k_1^2 - 24 j_0^3 d k_1^2 - 
    72 j_0^2 k_0 d k_1^2 - 72 j_0^2 j_1 d k_1^2\\
- 144 j_0 k_0 j_1 d k_1^2 - 144 j_0 j_1^2 d k_1^2 - 
    144 k_0 j_1^2 d k_1^2 - 240 j_1^3 d k_1^2 + 
    28 j_0^2 d^2 k_1^2 + 84 j_0 j_1 d^2 k_1^2 + 
    168 j_1^2 d^2 k_1^2\\
+ 9 j_0^3 k_1^3 + 
    18 j_0^2 k_0 k_1^3 + 36 j_0^2 j_1 k_1^3 + 
    48 j_0 k_0 j_1 k_1^3 + 90 j_0 j_1^2 k_1^3 + 
    60 k_0 j_1^2 k_1^3 + 180 j_1^3 k_1^3 - 24 j_0^2 d k_1^3\\
\left. - 96 j_0 j_1 d k_1^3 - 240 j_1^2 d k_1^3 + 6 j_0^2 k_1^4 + 
    30 j_0 j_1 k_1^4 + 90 j_1^2 k_1^4 \right)
\endmultline$$}

(V) Finally there are contributions from truncations, as explained in
Fact~5 of \S2
and Proposition~3.5.
A truncation determines two numbers $\ell$, $W$, and germs whose
limits $C_\sigma$ consist of unions of 4@-tangent conics and a
multiple of the distinguished tangent line; let $s_i$ denote
the multiplicities of the conics in $C_\sigma$, and write $S=\sum
s_i$. Then the contribution of the germ is
$$\ell W\left(192(S^7-\sum s_i^7)-288\,d(S^6-\sum
s_i^6)+112\,d^2(S^5-\sum s_i^5)\right)\quad.$$


\head \S4. Ordinary multiple points, and multiplicativity of adjusted
predegree polynomials\endhead

In this section we give an illustration of the results of \S3 by
obtaining explicit expressions for contributions accounting for
ordinary multiple points. We say that $p$ is an ordinary multiple
point for $C$ if $C$ has nonsingular branches with distinct tangent
directions at $p$; in particular, we allow branches to have flexes of
arbitrary order at $p$, or to be (reduced) lines. We also discuss to
what extent adjusted predegree polynomials are multiplicative with
respect to union of transversal curves.

\subheading{\S4.1. Ordinary multiple points} It is clear that ordinary
multiple points do not contribute components of type~V, since there is
only one branch along any direction of the tangent cone. The
contribution of an ordinary multiple point is therefore due to 1@-PS
germs, that is, components of type~III and IV.

\proclaim{Proposition~4.1} Let $p$ be an ordinary multiple point of
$C$, of multiplicity $m$, and for all lines $\ell$ tangent to a
non@-linear branch of $C$ at $p$ let $r_\ell$ be the intersection
multiplicity of $\ell$ and $C$ at $p$. Then the multiplicative
contribution to the adjusted predegree polynomial of $C$ due to $p$ is
given by{\eightpoint
$$\multline
\left(1-m^2 (m-1)(m-2)(m^2+3 m-3)\left(\frac{H^6}{720}-\frac{m
H^7}{840}+ \frac{m^2 H^8}{1920}\right)\right)\\
\cdot \prod_\ell 
\left(1-r_\ell (2 - 3 r_\ell + r_\ell^2 - 12 m + 3 r_\ell m + 6
m^2)\frac{H^6}{6!} + 
3 r_\ell (-12 + 2 r_\ell - 2 r_\ell^2 + r_\ell^3 \right.\\
+ 10 m - 8 r_\ell m +3 r_\ell^2 m - 20 m^2 + 6 r_\ell m^2 + 10
m^3)\frac{H^7}{7!} - 3 r_\ell (-64 + 2 r_\ell^2 - 3 r_\ell^3 + 2
r_\ell^4\\
\left. + 10 r_\ell m - 12 r_\ell^2 m + 6 r_\ell^3 m + 30 m^2 - 30
r_\ell m^2 + 12 r_\ell^2 m^2 - 60 m^3 + 20 r_\ell m^3 + 30
m^4)\frac{H^8}{8!}\right)
\endmultline$$}
where the $\prod$ is over all lines $\ell$ tangent to non@-linear
branches of $C$ at $p$.
\endproclaim
Note that linear branches do not appear directly in this formula,
although they have impact on the contribution by affecting $m$ and the
intersection multiplicities.

\demo{Proof} The first factor is the contribution of type~III, as in
Example~3.3.
According to Fact~4(ii) in \S2,
the other contributions from $p$ are due to the individual tangent
lines to the branches. Let $\ell$ be a line in the tangent cone to $C$
at $p$, and consider the branch of $C$ tangent to $\ell$ at $p$. We
note that

---if the branch is a line, $\ell$ does not contribute to the
   a.p.p.; indeed, the corresponding Newton polygon has no sides of
   slope strictly between $-1$ and $0$;

---if the branch is not a line, and has intersection multiplicity $k$
with $\ell$, then the corresponding Newton polygon has exactly one
side of slope strictly between $-1$ and $0$; this
side has vertices $(m-1,1)$ and $(r_\ell,0)$, where $r_\ell=m-1+k$ is
the intersection multiplicity of $\ell$ and $C$ at $p$. 

Applying Proposition~3.4
gives the contribution of type~IV due to $\ell$ in terms of $m$ and
$r_\ell$: this is the factor corresponding to $\ell$ in the
statement.\qed\enddemo

To state the result differently, let $e_i$ be the elementary symmetric
functions in the intersection multiplicities of $C$ with the tangent
lines to the non@-linear branches to $C$ at $p$. Then the multiplicative
contribution of $p$ to the a.p.p.~of $C$ is{\eightpoint
$$\multline
\bigg(1+(-2 e_1 + 3 e_1^2 - e_1^3 - 6 e_2 + 3 e_1 e_2 - 3 e_3 + 12 e_1
m - 3 e_1^2 m + 6 e_2 m + 6 m^2 - 6 e_1 m^2 - 15 m^3\\
+ 10 m^4 - m^6)\frac{H^6}{6!}
+(-36 e_1 + 6 e_1^2 - 6 e_1^3 + 3 e_1^4 - 12 e_2 + 18 e_1 e_2
- 12 e_1^2
e_2 + 6 e_2^2 - 18 e_3 + 12 e_1 e_3 - 12 e_4\\
+ 30 e_1 m - 24 e_1^2 m +
  9 e_1^3 m + 48 e_2 m - 27 e_1 e_2 m + 27 e_3 m - 60 e_1 m^2 + 18
e_1^2 m^2 - 36 e_2 m^2 - 36 m^3\\
+ 30 e_1 m^3 + 90 m^4 - 60 m^5 + 6
m^7) \frac{H^7}{7!}
+(192 e_1 - 6 e_1^3 + 9 e_1^4 - 6 e_1^5 + 18 e_1 e_2 - 36 e_1^2 e_2 + 30
e_1^3 e_2 + 18 e_2^2\\
- 30 e_1 e_2^2 - 18 e_3 + 36 e_1 e_3 - 30 e_1^2
e_3 + 30 e_2 e_3 - 36 e_4 + 30 e_1 e_4 - 30 e_5 - 30 e_1^2 m + 36
e_1^3 m - 18 e_1^4 m + 60 e_2 m\\
- 108 e_1 e_2 m + 72 e_1^2 e_2 m - 36
e_2^2 m + 108 e_3 m - 72 e_1 e_3 m + 72 e_4 m - 90 e_1 m^2 + 90 e_1^2
m^2 - 36 e_1^3 m^2\\
- 180 e_2 m^2 + 108 e_1 e_2 m^2 - 108 e_3 m^2 + 180
e_1 m^3 - 60 e_1^2 m^3 + 120 e_2 m^3 + 126 m^4 - 90 e_1 m^4 - 315 m^5\\
+ 210 m^6 - 21 m^8)\frac{H^8}{8!}\bigg)
\endmultline$$}

\example{Example~4.1} 
Suppose $p$ is an ordinary  {\it node\/} such that 
both branches of $C$ at $p$ intersect the respective tangent lines with
multiplicity exactly~2 at $p$. Then $p$ contributes
$$1-\frac{H^6}{6}+\frac{101 H^7}{280}-\frac{25 H^8}{64}$$
to the a.p.p.~(set $m=2$, $e_1=3+3$, $e_2=3\cdot 3$, $e_3=e_4=e_5=0$
in the previous formula).
Since $p$ `absorbs' 6 ordinary inflection points, the
predegree polynomial for a curve of degree $d\ge 3$ with
$n$ such nodes and only ordinary flexes is{\eightpoint
$$\multline
\exp(d H)\cdot\left(1-2d\left(\frac{H^5}{20}-\frac{(5d+18)H^6}{360}
+\frac{(9d+8) H^7}{420}-\frac{d H^8}{60}\right)\right)\\
\cdot\left(1-\frac{H^6}{42}+\frac{3 H^7}{70}-\frac{197
H^8}{4480}\right)^{3 d(d-2)-6n} \cdot \left(1-\frac{H^6}{6}+\frac{101
H^7}{280}-\frac{25 H^8}{64}\right)^n
\endmultline$$}
(The term following the exponential is the contribution as in
Example~3.2;
the next term accounts for the flexes, obtained by setting $k=3$ in
Example~3.4.)
The predegree of such a curve is therefore 
$$d^8-1372 d^4+7992 d^3-15879 d^2+10638 d-24 n \left(35 d^2-174 d+213
\right)\quad.$$
For instance, the degree of the orbit closure of 
a quartic of this kind is $14280-1848 n$;
the predegree of the orbit closure of a rational plane curve
of this kind is
$$d^8-1792 d^4+11340 d^3-25539 d^2+22482 d-5112\quad.$$
\endexample

\example{Example~4.2} An ordinary multiple point $p$ of multiplicity $m$,
and such that each branch is smooth, non@-linear, and does not have an
inflection point at $p$ contributes{\eightpoint
$$\multline
\left(1-\frac{m(m^3+m^2+m+16)H^6}{6!} +\frac{3(2 m^5+2 m^4+2 m^3+37
m^2+16 m+11)H^7}{7!}\right.\\
\left.-\frac{21(m^6+m^5+m^4+21 m^3+13 m^2+17 m+9)H^8}{8!}
\right)^{m(m-1)}\quad.
\endmultline$$}
Using that such a point absorbs $3 m (m-1)$ flexes,
one then sees that the contribution to the predegree of a
curve of degree $d$ due to such a point is{\eightpoint
$$\multline
-m(m-1)(21 m^6-48 d m^5+21 m^5+28 d^2 m^4-48 d m^4+21 m^4+28 d^2
m^3-48 d m^3\\
+441 m^3+28 d^2 m^2-888 d m^2+273 m^2+448 d^2 m-384 d m+357 m
-1260 d^2+4920 d-5130). 
\endmultline$$}
For instance, a general quartic curve with a triple point has
predegree $14280-3\cdot2\cdot 1890=2940$.
\endexample

\example{Example~4.3} A {\it biflecnode\/} is an ordinary node at which
both branches have an ordinary inflection point; its contribution is
$$1-\frac{H^6}3+\frac{88 H^7}{105}-\frac{15 H^8}{14}$$
(set $m=2$, $e_1=4+4$, $e_2=4\cdot 4$, $e_3=e_4=e_5=0$ in the formula
given above). Using that such a point absorbs 8~flexes, we get that a
biflecnode corrects the predegree for a curve of degree~$d$ by
$$-24(140 d^2-832 d+1209)\quad.$$
For instance, the quartic with equation
$$x^2 y^2 + x^2 z^2 + y^2 z^2=0$$
has three biflecnodes and 24 automorphisms, hence its orbit closure
has predegree $\frac{14280-3\cdot 2904}{24}=232$. 
As it happens, this orbit closure is isomorphic to the moduli space
of semistable vector bundles on $\P^2$ of rank 2 with Chern classes
$c_1=-1$ and $c_2=3$, as Hulek proved \cite{H}.
It follows that the corresponding Donaldson invariant of $\P^2$
equals 232, in agreement with \cite{K-L}.
\endexample

\example{Example~4.4} Suppose $p$ is an ordinary node for which one branch
is a {\it line,\/} and the other intersects its tangent line with
multiplicity $k$ at $p$. Then $p$ contributes{\eightpoint
$$\multline
1-\frac{(k+1)(k+2)(k+3) H^6}{6!}
+ \frac{3(k+1)(k^3+7k^2+21k+23)H^7}{7!}\\
-\frac{3(k+1)(k+3)^2(2k^2+5k+17)H^8}{8!}
\endmultline$$}
(use $m=2$, $e_1=k+1$, $e_2=e_3=e_4=e_5=0$ in the formula given above).
For $k=2$, the contribution is
$$1-\frac{H^6}{12}+\frac{101 H^7}{560}-\frac{25 H^8}{128}\quad;$$
of course this is the square root (modulo $H^9$) of the contribution
for a node given in Example~4.1.
\endexample

\subheading{\S4.2. Multiplicativity of adjusted predegree polynomials}
It is natural to ask whether the predegree information behaves well
with respect to unions of curves. This is another advantage of
adjusted predegree polynomials over other ways to assemble this
enumerative information: adjusted predegree polynomials are
multiplicative under unions of curves, up to correction terms
independent of the degree(!),
accounting for the ways in which the curves meet.
No such structure is visible at the level of degrees or predegrees
alone.

As a representative example, we let $C_1$, $C_2$ be arbitrary {\it
reduced\/} curves, meeting {\it transversally\/} at nonsingular
points, and we further assume that such points are not inflection
points for either curve.
Let $C'_i$ resp.~$L_i$ be the union of the non@-linear
resp.~linear components of $C_i$. Let $I=\#(C'_1\cap C'_2)$,
$J=\#((C'_1\cap L_2)\cup (C'_2\cap L_1))$.

\proclaim{Proposition~4.2} Let $P_{C_1}(H)$, $P_{C_2}(H)$ be the adjusted
predegree polynomials of $C_1$, $C_2$. Then the adjusted predegree
polynomial of their union $C=C_1\cup C_2$ is
$$\multline 
P_C(H)=P_{C_1}(H)\cdot P_{C_2}(H)\\
\cdot\left(1-\frac{H^6}{9}+\frac{11H^7}{40}-\frac{311 H^8}{960}\right)^I
\cdot\left(1-\frac{H^6}{24}+\frac{7H^7}{60}-\frac{13 H^8}{80}\right)^J
\quad.
\endmultline$$
\endproclaim

\demo{Proof} The main remark is that the components of the projective
normal cone for $C_1\cup C_2$ arise from features of $C_1$, $C_2$ and
from the points of intersection of the two curves; an analysis of the
components leads to the formula of the statement. We go through this
analysis here as a template for similar computations. 

As pointed out in Example~2.3,
the intersection of two lines does not contribute components. Using the
formulas given in Examples~4.1--4.4
to evaluate the contribution of the transversal intersections of two
curves at non@-flex point, and of a line and a curve at a non@-flex
point, we can write
$$\multline
P_C(H)=\exp((d_1+d_2) H)\left(1+L_{C'_1}(C)+L_{L_1}(C) + L_{C'_2}(C)
+ L_{L_2}(C)\right)\\
\cdot \left(1+L_{\text{local}}(C_1)\right)\left(1+
L_{\text{local}}(C_2)\right)\left(1-\frac{H^6}{12}+\frac{101
H^7}{560}- \frac{25 H^8}{128}\right)^{2I+J}
\endmultline$$
where $d_i=\deg C_i$, and $L_{\dots}$ denote the various correction
terms, with hopefully evident notations: for example,
$L_{\text{local}}(C_1)$ stands for the term arising from all local
features of $C_1$. It is crucial here to recall (cf.~Lemma~3.3.1)
that such local terms do not depend on other features of the curve; so
the contribution of a local term is the same whether viewed in $C_i$
or in $C$. (This is not the case for `global' terms!) With the same
notations we can write
$$P_{C_i}(H)=\exp(d_i H) \left(1+L_{C'_i}(C_i)+L_{L_i}(C_i)\right)
\left(1+L_{\text{local}}(C_i)\right)$$
and therefore the ratio $\frac{P_C(H)}{P_{C_1}(H) P_{C_2}(H)}$ is
expressed by{\eightpoint
$$\frac{(1+L_{C'_1}(C)+L_{L_1}(C) + L_{C'_2}(C)+
L_{L_2}(C))}{(1+L_{C'_1}(C_1)+L_{L_1}(C_1))
(1+L_{C'_2}(C_2)+L_{L_2}(C_2))} \left(1-\frac{H^6}{12}+\frac{101
H^7}{560}- \frac{25 H^8}{128}\right)^{2I+J}\quad.$$}
Lastly, we note that in evaluating this term we may assume that 
each line meets the rest of $C$ transversally at non-inflection points:
indeed, the terms arising from
special positions of the lines can be evaluated locally, so they can
be incorporated in the $L_{\text{local}}$ terms. All the terms in this
expression can then be evaluated very simply by Propositions~3.1 and
3.2,
giving the stated result.\qed\enddemo

\example{Example~4.5} If both $C_1$, $C_2$ are unions of lines, then
multiplicativity holds `on the nose', since $I=J=0$ in that case. This
in fact holds for non@-reduced configurations of lines as well,
cf.~Corollary~2.7 in \cite{A-F4}.
\endexample

\example{Example~4.6} The union of a general curve $C$ of degree $d\ge 2$ and
a general transversal line has predegree polynomial
$$P_C(H)\cdot\left(1+H+\frac{H^2}2\right)\cdot
\left(1-\frac{H^6}{24}+\frac{7H^7}{60}-\frac{13 H^8}{80}\right)^d
\quad,$$
where $P_C(H)$ is the predegree polynomial of a general curve (computed
in Example~1.1).
For $d=2$, this yields
$$1+3H+\frac{9H^2}2+\frac{13 H^3}3+3 H^4+\frac{7 H^5}5+\frac{19
H^6}{60} +\frac{H^7}{60}\quad,$$
detecting that the union of a conic and a transversal line has orbit
closure of dimension~7 and degree $\frac{7!}{60\cdot 4}=21$.
This agrees of course with the na\"\i ve combinatorial count, since the
orbit of the union of a conic and a transversal line is in fact the
set of all such curves; the degree is then the number of curves
through 7 general points, that is $\binom 72=21$ (the line must
contain two of the points, and the conic is then determined by the
other five). 

Combinatorics would not suffice to compute e.g.~the degree for the
union of a general {\it cubic\/} and a general transversal line;
according to the formula given above, this is~$8568$.
Note that these computations do depend on whether the intersection
points are or are not inflection points for the branches. Using the
formula given in Example~4.4,
one obtains that the predegree of the union of a general cubic and
a general transversal line {\it through a flex\/} of the cubic is
$8040$.
\endexample

\example{Example~4.7} The union of two transversal conics has a.p.p.~given
by
$$\multline
\left(1+2 H+2 H^2+\frac{4 H^3}3+\frac{2
H^4}3+\frac{H^5}{15}\right)^2 \cdot \left(1-\frac{H^6}{9}
+\frac{11H^7}{40}-\frac{311 H^8}{960}\right)^4\\
=1+4 H+8 H^2+\frac{32 H^3}3+\frac{32 H^4}3+\frac{122 H^5}{15}
+\frac{64 H^6}{15}+\frac{41 H^7}{30}+\frac{41 H^8}{240}\quad,
\endmultline$$
hence predegree $6888$.
\endexample

The reader will have no difficulties adapting the argument in the
proof of Proposition~4.2
to compute terms accounting for other kinds of intersections. For
example, a point of simple tangency of a line with a curve gives a
correction term
$$1-\frac{H^6}6+\frac{7 H^7}{15}-\frac{13 H^8}{20}$$
to the polynomial of the union of the curve and the line (note that
this is the 4@-th power (modulo $H^9$) of the contribution for a point
of transversal intersection of a line with a curve. We don't have a
conceptual explanation for this phenomenon).
Thus, the predegree polynomial for the union of a smooth conic and a
tangent line is{\eightpoint
$$\multline
\left(1+2 H+2 H^2+\frac{4 H^3}3+\frac{2 H^4}3+\frac{H^5}{15}\right)
\cdot \left(1+H+\frac{H^2}2\right) \cdot \left(1-\frac{H^6}6+\frac{7
H^7}{15}-\frac{13 H^8}{20}\right)\\
=1+3 H+\frac{9 H^2}2+\frac{13 H^3}3+3 H^4+\frac{7 H^5}5+\frac{7
H^6}{30}\quad:
\endmultline$$}
the orbit closure has dimension~6 and degree $\frac{6! 7}{30 \cdot
4}=42$, as expected. 


\head \S5. Irreducible singularities\endhead

Our last and most substantial example illustrating the algorithm
implicitly described in \S\S2@-3
will be the computation of the contribution to the adjusted predegree
polynomial due to an arbitrary  irreducible singularity $p$
on a curve~$C$.

It is well@-known
that $C$ can be described at such a point by its {\it Puiseux
expansion\/}
$$\left\{\aligned
z&=(a_n t^n+\dots+)\,\,a_{e_1} t^{e_1}+\dots+a_{e_r} t^{e_r}\\
y&=t^m
\endaligned\right.$$
where: $m=$ the multiplicity of $C$ at $p$; $n=$ the intersection
multiplicity of $C$ and the tangent line $z=0$ at $p$; all
exponents are positive integers, and $m<n\le e_1<\dots< e_r$; and the
coefficients $a_{e_i}$ of the `essential' terms are nonzero. An
exponent (or the corresponding term in the expansion) is `essential'
if it is {\it not\/} a multiple of the greatest common divisor of $m$
and the exponents preceding it; the $()$ in the expansion collects all
non@-essential terms. The term $a_n t^n$ will be essential 
if and only if $n$ is not a
multiple of $m$; note that $e_1=n$ in that case.

We also need the numbers
$$d_i=gcd(m,e_1,\dots,e_i)\quad;$$ 
thus $d_0=m$, and $d_r=1$. Note that we allow for the possibility
$m=1$, $r=0$; that is, there may be {\it no\/} essential terms in the
expansion.

We will see that the contribution of $p$ to the a.p.p.~for $C$ {\it
depends only on $m$, $n$, and the essential exponents $e_1$, \dots,
$e_r$.\/}

An alternative terminology to describe the same information is that of
{\it Puiseux pairs:\/} the singularity is described by the pair
$(m,n)$, and by $r$ Puiseux pairs $(m_1,n_1)$, \dots, $(m_r,n_r)$,
where
$$\left\{\aligned
d_i&=m_{i+1}\cdots m_r\\
e_i&=n_i d_i
\endaligned\right.\quad.$$
Thus for example a non@-singular inflection point of order $k$ is
described by
$$(1,k)$$
and has no Puiseux pairs ($r=0$, no essential exponents,
$d_0=1=m$); an ordinary cusp $y^n=z^m$ ($m$, $n$ coprime) is described
by
$$(m,n)\quad;\quad (m,n)$$
and has {\it one\/} Puiseux pair ($r=1$, $e_1=n$, $d_0=m$,
$d_1=1$). The formula given below
implies that the correction due to $p$ only depends on $m$, $n$, and
the Puiseux pairs of $C$ at $p$.

This result is most easily stated in terms of the numbers $d_i$,
$e_i$. We let 
$$P(a,b)=\frac{a^2 b^2}{(1+a k)^3(1+b
k)^3}-\frac{4}{(1+k)^3(1+2k)^3}\quad,$$
where $k$ is an indeterminate, and set $e_0=n$, $e_{r+1}=0$ for
convenience.

\proclaim{Theorem~5.1} With notations as above, the contribution of $p$ to
the adjusted predegree polynomial of $C$ is
$$1-\left\{\bigg(mnP(m,n)+\sum_{j=0}^{r}(e_{j+1}-e_j) d_j P(d_j,2 d_j)
\bigg)\cdot \left(\frac{k^2 H^6}{6!}+\frac{k H^7}{7!}+\frac{H^8}{8!}
\right)\right\}_2$$
where $\{\}_2$ denotes the coefficient of $k^2$ in the expansion of
the term within $\{\}$.\endproclaim

Before proving this formula, we illustrate it with a few explicit
examples. For these we will need the number of flexes absorbed
by the singularity; remarkably, this number can be expressed by a
formula somewhat analogous to the one given in Theorem~5.1:
$$(3mn-2m-2n)+3\sum_{j=0}^r(e_{j+1}-e_j)(d_j-1)$$
(cf.~\cite{B-K}, \S9.1, Thm.~2 and \cite{Oka}, \S2). 
The correction term that would be due to the flexes absorbed by $p$ if
$p$ were not present is, according to Theorem~5.1,
{\eightpoint
$$1-\left\{\bigg((3mn-2m-2n)+3\sum_{j=0}^r(e_{j+1}-e_j)(d_j-1)\bigg)
3P(1,3)\cdot \left(\frac{k^2 H^6}{6!}+\frac{k H^7}{7!}+\frac{H^8}{8!}
\right)\right\}_2\quad.$$}

\example{Example~5.1} A nonsingular point has no Puiseux pairs, and
$(m,n)=(1,k)$, where $k=$ the order of contact with the tangent
line. By Theorem~5.1,
its contribution is{\eightpoint
$$\multline
1-\left\{kP(1,k)\cdot \left(\frac{k^2 H^6}{6!}+\frac{k
H^7}{7!}+\frac{H^8}{8!}\right)\right\}_2\\
=1-\frac{k(k-2)(k+2)H^6}{720}+\frac{k(k-2)(k^2+3k+6)H^7}{1680}
-\frac{k(k-2)(2k^3+7k^2+16k+32)H^8}{13440}
\endmultline$$} 
in agreement with Example~3.4.
Note that this contribution is automatically trivial if $k=2$, that
is if the point is not an inflection point for $C$.
\endexample

Assume next that $p$ has exactly one Puiseux pair $(m_1,n_1)$. With
notations as above, necessarily $m_1=m$; and $d_0=m$, $d_1=1$;
$e_0=n$, $e_1=n_1$, $e_2=0$. According to Theorem~5.1,
the contribution of $p$ is{\eightpoint
$$\multline
1-\frac{m(4 m^4 (n_1-n) + m^2 n^3 - 4 n_1)H^6}{6!}
+\frac{3 m (12 m^5 (n_1-n) + m^3 n^3 + m^2 n^4 - 12 n_1)H^7}{7!}\\
-\frac{3 m (64 m^6 (n_1-n) + 2 m^4 n^3 + 3 m^3 n^4 + 2 m^2 n^5 - 64
n_1)H^8}{8!}
\endmultline$$}
\example{Example~5.2} For an ordinary $(m,n)$@-cusp 
(see above) we find{\eightpoint
$$1-\frac{m n(m^2 n^2-4)H^6}{6!}+\frac{3m n(m^3 n^2+m^2 n^3-12)
H^7}{7!} -\frac{3 m n(2 m^4 n^2+3 m^3 n^3+2 m^2 n^4-64)H^8}{8!}$$}
For instance, an ordinary $(2,3)$ cusp contributes
$$1-\frac{4 H^6}{15}+\frac{3 H^7}5-\frac{19 H^8}{28}\quad;$$
using that such a cusp absorbs 8 flexes, we get that an ordinary cusp
corrects the predegree of a curve of degree $d\ge3$ by
$$-72(28 d^2-144 d+183)\quad.$$
Thus a generic cuspidal quartic has predegree $14280-3960=10320$,
etc. Note that for a cuspidal {\it cubic\/} this gives a `predegree'
of $216-216=0$; this is because cuspidal cubics have small
orbits. According to the formulas given above, the a.p.p.~of a
cuspidal cubic is
$$1+3 H+\frac{9 H^2}2+\frac{9 H^3}2+\frac{27 H^4}8+\frac{69 H^5}{40}
+\frac{3 H^6}8+\frac{H^7}{70}\quad,$$
yielding a degree of $\frac{7!}{70\cdot 3}=24$, as expected.
\endexample

\example{Example~5.3} {\it Characteristic numbers.\/} An enumerative
problem that has received a good deal of attention both in the 19th
century and in the recent past is that of computing the {\it
characteristic numbers\/} of various families of plane curves, that
is, the number of curves belonging to the family, containing a
collection of general points, and tangent to a collection of general
lines. This problem is in general surprisingly challenging, even for
curves of very low degree.

We note here that the top characteristic number of the (family of
curves parameterized by the) orbit closure of $C$ is the degree of the
orbit closure of the dual curve $C^\vee$: hence the results of this
paper allow us in principle to compute the `top' characteristic number
of the orbit closure of an arbitrary curve, i.e., the number of
translates of the curve which are tangent to a maximal number of
general lines.

For example, consider the orbit closure of a nonsingular cubic curve
$C$, that is, the closure of the set of cubic curves with a given
$j$@-invariant. Its top characteristic number is the degree of
the orbit closure of a sextic with 9~cusps; now Example~5.2
lets us compute the predegree of this orbit closure:
$$\multline
\text{predegree of a general sextic} - \text{contributions from 9
cusps}\\
= 1119960-9\cdot 23544 =908064\quad.
\endmultline$$
For $j\ne 0,1728$, the stabilizer of $C$ consists of 18 elements; 
thus there are $\frac{908064}{18}=\oldnos{50,448}$ cubics with fixed
$j$ invariant $\ne 0,1728$ and tangent to 8~lines in general position.
For $j=0$, resp.~$j=1728$ the extra automorphisms of $C$ correct this
number to $\frac{50448}{3}=\oldnos{16,816}$,
$\frac{50448}{2}=\oldnos{25,224}$, respectively. These results agree
with the more direct computations in \cite{A}.

Similarly, the number of nodal cubics tangent to $8$ lines in general
position is the degree of the orbit closure of the dual of a
nodal cubic, that is, a quartic with three cusps:
$$\frac{14280-3\cdot 3960}6=\oldnos{400}\quad.$$
Of course this also agrees with the classical result (cf.~for example
\cite{S}).

It is curious to observe that the dual of a nodal cubic can also be
interpreted as a sextic consisting of a quartic with three cusps and a
double bitangent line, in the sense that this is what the dual of a
nonsingular cubic $C$ degenerates to as $C$ degenerates to a nodal
cubic. Arguing as in \S4 to account for the contribution of the
double line, we compute that the predegree of the orbit closure of
such a sextic is $302668$; as the stabilizer of a nodal cubic has~6
elements,
this gives $\oldnos{50,448}$ as the `top characteristic number' of a
nodal cubic. This number counts the 400 curves tangent to 8~lines as
well as contributions from curves whose node is on one of the lines;
the fact that this number agrees with the characteristic number for
cubics with $j<\infty$ was already observed in \cite{A}, end of \S3.

Apart from these and a few other instances (for example conics, or
cuspidal cubics), the characteristic numbers that can be obtained by
applying the results in this paper are, to our knowledge, new. For
example, so is the number \oldnos{406,758,744} of nonsingular {\it
quartics\/} with fixed general modulus
and tangent to 8~lines in general position.
\endexample

\example{Example~5.4} The quartic curves 
$$(y^2-x z)^2=y^3 z\quad; \quad (y^2-x z)^2=y z^3$$
have a singularity at $(1:0:0)$ described by $(m,n)=(2,4)$, and
Puiseux pair $(2,k)$ for $k=5$, $7$ respectively. Using the formula
given above, and that these points absorb $3k$ flexes, we find that
these singularities correct the predegree of the quartics on which
they lie by $-1785 k$.

These singularities are analytically isormorphic to $z^2=y^k$ 
(cf.~Example~3.5).
Remarkably, the same correction term applies for quartics with a point
analytically isomorphic to $z^2=y^k$ also in the non@-irreducible
cases $k=4$, $6$, $8$ (as may be computed explicitly using
Propositions~3.4 and 3.5).
For $k=8$ the corresponding quartic is $(y^2-x z)^2=z^4$, that is the
union of two quadritangent conics (cf.~\S4.1 in \cite{A-F3});
the formula gives $14280-1785\cdot 8=0$, as expected since unions of
quadritangent conics have small orbits.

The case $k=4$ can also be analyzed by the same method, and gives a
correction of $-1785\cdot 4=-7140$. Thus a general tacnodal quartic has
predegree $14280-7140=7140$, that is, precisely half of the predegree
of a general quartic. This latter fact can also be explained
conceptually by studying the behavior of the predegree along
families of curves, but we will not pursue this approach here.
\endexample

\demo{Proof of Theorem~5.1}
The formula given in the theorem is obtained by evaluating explicitly
the contributions of type~IV and V, using Proposition~3.4 and 3.5.
The main subtlety lies in the fact that both these contributions are
affected by whether $n$ is an essential exponent or not; as we will
see, the amounts by which they are affected precisely compensate
each other, so that both cases lead to the same formula. 

We consider contributions of type~IV first. If $d'=gcd(m,n)$, and
$m'=m/d'$, $n'=n/d'$, then the only 1@-PS germ giving a contribution
is
$$\pmatrix
1 & 0 & 0\\
0 & t^{n'} & 0\\
0 & 0 & t^{m'}
\endpmatrix$$
yielding a limit
$$\left(y^{n'}-*x^{n'-m'} z^{m'}\right)^{d'} x^{d-n}=0\quad,$$
corresponding to the side in the Newton polygon joining vertices
$(0,m)$ and $(n,0)$. Using Proposition~3.4,
this gives a contribution of{\eightpoint
$$\multline
1-mn\left(\frac{(m^2 n^2-4 {d'}^4)H^6}{6!}-\frac{3(m^3 n^2+m^2 n^3-12
{d'}^5) H^7}{7!}\right.\\
\left.+\frac{3(2 m^4 n^2+3 m^3 n^3+2 m^2 n^4-64
{d'}^6)H^8}{8!}\right)
\endmultline$$}
which is checked to equal{\eightpoint
$$\multline
1-\left\{mn\big(P(m,n)-P(m,2m)\big)\cdot \left(\frac{k^2 H^6}{6!}+
\frac{k H^7}{7!}+\frac{H^8}{8!} \right)\right\}_2\\
-mn \left(\frac{(m^4-{d'}^4) H^6} {180}-\frac{(m^5-{d'}^5)H^7}
{140}+\frac{(m^6-{d'}^6)H^8} {210}\right)\quad.
\endmultline$$}
Here $d'=m$ if $n$ is a multiple of $m$ (in which case the last
summand vanishes), while $d'=d_1=gcd(m,e_1)$ if $n=e_1$ is essential. 

Moving on to the component of type~V, the data describing the
singularity determines
the structure of the formal branches of the curve at
$p$. Schematically, here is how they group:
$$\epsffile{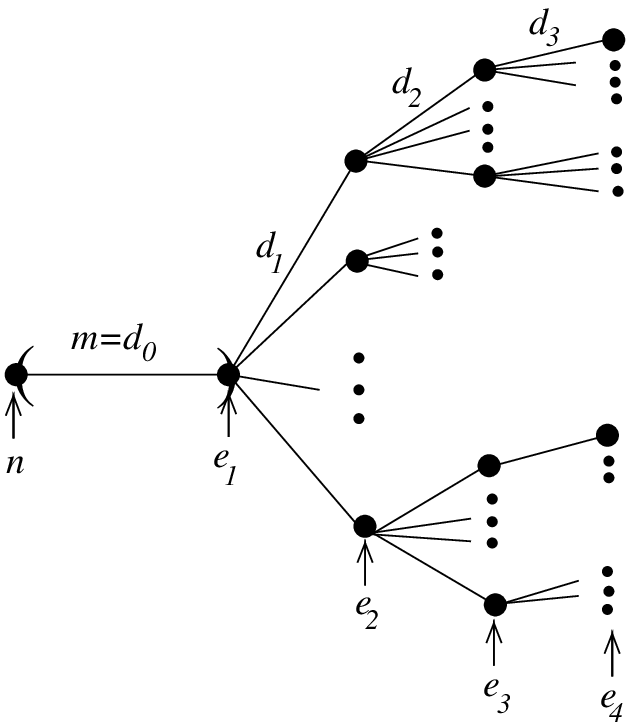}$$
If $n$ is not essential, $m=d_0$ branches will run parallel from the
beginning of the expansion up to the first essential exponent $e_1$;
if $n$ is essential, the branching starts immediately at $n=e_1$.
In both cases, at $e_1$ the branches divide into $d_0/d_1$ groups of
$d_1$ parallel branches each; at $e_2$, each set of $d_1$ branches
splits into $d_1/d_2$ groups of $d_2$ parallel branches, and so on. At
the last essential exponent $e_r$, the splitting produces $m$ distinct
simple branches.

This gives us the data needed to apply Proposition~3.5.
Note that $e_1$ yields a `truncation' in the sense of Fact~5 of \S2
{\it only if\/} $n$ is not an essential exponent: if $n=e_1$ is
essential, then the expansion starts at $e_1$ and in particular $e_1$
is not greater than the first exponent.
If $n$ is not essential, the truncation at $e_1$ contributes in the
terminology of Proposition~3.5
a term with $\ell=1$, $W=e_1$, $S=m$, and $s_i=d_1$, giving
$$-me_1 \left(\frac{(m^4-d_1^4) H^6} {180}-\frac{(m^5-d_1^5)H^7}
{140}+\frac{(m^6-d_1^6)H^8} {210}\right)\quad;$$
if $n$ is essential, there is no such contribution. Adding this to the
contribution of type~IV computed above, we obtain {\it in both
cases\/}
$$1-\left\{mn\big(P(m,n)-P(m,2m)\big)\cdot \left(\frac{k^2 H^6}{6!}+
\frac{k H^7}{7!}+\frac{H^8}{8!} \right)\right\}_2-K_1$$
where
$$K_1=me_1 \left(\frac{(m^4-d_1^4) H^6} {180}-\frac{(m^5-d_1^5)H^7}
{140}+\frac{(m^6-d_1^6)H^8} {210}\right)\quad.$$

The contribution due to truncation at $e_j$, $j\ge 2$, is given by
Proposition~3.5,
setting $\ell=\frac m{d_{j-1}}$ (the least integer such that
$\ell\frac{e_1}m$, \dots, $\ell\frac{e_{j-1}}m$ are integers),
$$W=\sum_{k=1}^{j-1} (d_{k-1}-d_k)\frac{e_k}m+d_{j-1}\frac{e_j}m$$
(keeping track of the exponents at which formal branches start
differing),
and $S=d_{j-1}$, $s_i=d_j$. If $K_j$ denotes this (additive)
contribution, one checks by induction that if there are $r$ Puiseux
pairs (so that $d_r=1$)
$$\sum_{j=2}^rK_j=\left\{\sum_{j=2}^r e_j\big(d_{j-1}P(d_{j-1},2
d_{j-1})-d_j P(d_j,2 d_j)\big)\left(\frac{k^2 H^6}{6!}+ \frac{k
H^7}{7!}+\frac{H^8}{8!}\right)\right\}_2$$
(note: this equality does not hold if $d_r$ is not assumed to
equal~1!).
The whole contribution is therefore given by
$$1-\left\{mn\big(P(m,n)-P(m,2m)\big)\cdot \left(\frac{k^2 H^6}{6!}+
\frac{k H^7}{7!}+\frac{H^8}{8!}\right)\right\}_2
-\sum_{j=1}^rK_j\quad,$$
and the formula given in the statement is obtained by rearranging this
sum.\qed\enddemo

Formulas for {\it reducible\/} singularities can be obtained by using
Propositions~3.3, 3.4, and 3.5.
Unfortunately,
we haven't been able to find a simple statement in the style of
Theorem~5.1
and encompassing the most general case.

As a final comment, we note that a formula in the style of
Theorem~5.1
can be concocted to account for some `global' terms as well. For
example, the predegree of the orbit closure of a reduced curve of
degree $d$ and (for simplicity) including only points `of type
$(t^m,t^n)$' (that is, points described by the pair $(m,n)$ as above,
with no Puiseux pairs) is in fact given by
$$\multline
d^8-\bigg\{(1+dk)^8\bigg[\frac{4 d^2}{(1+k)^3(1+2 k)^3}\\
+\sum_{p\in C\,\text {of type }(t^m,t^n)}
mn\left(\frac{m^2n^2}{(1+m k)^3(1+nk)^3}- \frac{4}{(1+k)^3
(1+2 k)^3}\right)\bigg]\bigg\}_2\quad,
\endmultline$$
provided that the orbit closure has dimension~8. This formula should
be compared with the formula for the predegree of the orbit closure of
a {\it $d$@-tuple of points in $\P^1$\/} (cf.~\cite{A-F1}),
which can be written
$$d^3-\bigg\{(1+d k)^3\bigg[ \frac d{(1+k)^2} + \sum_{p\in C\,\text
{of type }(t^m)} m\left( \frac m{(1+mk)^2} - \frac
1{(1+k)^2}\right)\bigg]\bigg\}_1$$
(if the orbit closure has dimension~3), where a point `of type
$(t^m)$' is simply a point of multiplicity $m$ in the $d$@-tuple.

It is tempting to view these two formulas as shadows of a very
general, but as yet mysterious, theorem on degrees of orbit closures
of hypersurfaces in projective space.


\leftheadtext{References}
\rightheadtext{References}

\enddocument